\documentclass[10pt]{article} 
\usepackage{latexsym} 
\usepackage{amsmath}
\usepackage{amssymb} 
\usepackage{bm}
\usepackage{amscd} 
\usepackage{graphicx}
\begin{document} 
%%%%%%%%%%%%%%%%%%%%%%%%%%%%%%%%% 
%%%%%%%%%%%%%%%%%%%%%%%%%%%%%%%%%
\newtheorem{Th}{Theorem}[section]
\newtheorem{Cor}{Corollary}[section]
\newtheorem{Prop}{Proposition}[section]
\newtheorem{Lem}{Lemma}[section]
\newtheorem{Def}{Definition}[section]
\newtheorem{Rem}{Remark}[section]
\newtheorem{Ex}{Example}[section]
\newtheorem{stw}{Proposition}[section]

%Definitions of bet ent

\newcommand{\bet}{\begin{Th}}
\newcommand{\ent}{\stepcounter{Cor}
   \stepcounter{Prop}\stepcounter{Lem}\stepcounter{Def}
   \stepcounter{Rem}\stepcounter{Ex}\end{Th}}

%Definitions of bec enc bep enp bel enl
%bef enf ber enr bee ene

\newcommand{\bec}{\begin{Cor}}
\newcommand{\enc}{\stepcounter{Th}
   \stepcounter{Prop}\stepcounter{Lem}\stepcounter{Def}
   \stepcounter{Rem}\stepcounter{Ex}\end{Cor}}
\newcommand{\bep}{\begin{Prop}}
\newcommand{\enp}{\stepcounter{Th}
   \stepcounter{Cor}\stepcounter{Lem}\stepcounter{Def}
   \stepcounter{Rem}\stepcounter{Ex}\end{Prop}}
\newcommand{\bel}{\begin{Lem}}
\newcommand{\enl}{\stepcounter{Th}
   \stepcounter{Cor}\stepcounter{Prop}\stepcounter{Def}
   \stepcounter{Rem}\stepcounter{Ex}\end{Lem}}
\newcommand{\bef}{\begin{Def}}
\newcommand{\enf}{\stepcounter{Th}
   \stepcounter{Cor}\stepcounter{Prop}\stepcounter{Lem}
   \stepcounter{Rem}\stepcounter{Ex}\end{Def}}
\newcommand{\ber}{\begin{Rem}}
\newcommand{\enr}{
   %\stepcounter{Rem} 
   \stepcounter{Th}\stepcounter{Cor}\stepcounter{Prop}
   \stepcounter{Lem}\stepcounter{Def}\stepcounter{Ex}\end{Rem}}
\newcommand{\bee}{\begin{Ex}}
\newcommand{\ene}{
 %\stepcounter{Ex}
   \stepcounter{Th}\stepcounter{Cor}\stepcounter{Prop}
   \stepcounter{Lem}\stepcounter{Def}\stepcounter{Rem}\end{Ex}}
\newcommand{\Proof}{\noindent{\it Proof\,}:\ }

%%%%%%%%%%%%%%%%%%%%%%%%%%%%%%%%%%%%%%%%%
%Beginning of Local Definition
%Local definitions
\newcommand{\R}{{\mathbb{R}}}
\newcommand{\C}{{\mathbb {C}}}
\newcommand{\K}{{\mathbb {K}}}

\newcommand{\ZZ}{{\mathbb {Z}}}
\newcommand{\NN}{{\mathbb {N}}}
\newcommand{\M}{{\mathcal{M}}}
\newcommand{\KK}{{\mathbb {K}}}
\newcommand{\PP}{{\mathbb {P}}}
\newcommand{\OOO}{{\mathbb {O}}}
\newcommand{\xx}{{\mathbold{x}}}
\newcommand{\EE}{\mathbb{E}}
\newcommand{\QQ}{\mathbb{Q}}
\newcommand{\uuu}{\boldsymbol{u}}
\newcommand{\xxx}{\boldsymbol{x}}
\newcommand{\aaa}{\boldsymbol{a}}
\newcommand{\bbb}{\boldsymbol{b}}
\newcommand{\AAA}{\mbox{\bf{A}}}
\newcommand{\BBB}{\mbox{\bf{B}}}
\newcommand{\ccc}{\boldsymbol{c}}
\newcommand{\iii}{\boldsymbol{i}}
\newcommand{\jjj}{\boldsymbol{j}}
\newcommand{\kkk}{\boldsymbol{k}}
\newcommand{\rrr}{\boldsymbol{r}}
\newcommand{\FFF}{\boldsymbol{F}}
\newcommand{\yyy}{\boldsymbol{y}}
\newcommand{\ppp}{\boldsymbol{p}}
\newcommand{\qqq}{\boldsymbol{q}}
\newcommand{\nnn}{\boldsymbol{n}}
\newcommand{\vvv}{\boldsymbol{v}}
\newcommand{\eee}{\boldsymbol{e}}
\newcommand{\fff}{\boldsymbol{f}}
\newcommand{\www}{\boldsymbol{w}}
\newcommand{\0}{\boldsymbol{0}}
\newcommand{\lon}{\longrightarrow}
\newcommand{\ga}{\gamma}
\newcommand{\pa}{\partial}
\newcommand{\QED}{\hfill $\Box$}
\newcommand{\id}{{\mbox {\rm id}}}
\newcommand{\Ker}{{\mbox {\rm Ker}}}
\newcommand{\grad}{{\mbox {\rm grad}}}
\newcommand{\ind}{{\mbox {\rm ind}}}
\newcommand{\rot}{{\mbox {\rm rot}}}
\newcommand{\diver}{{\mbox {\rm div}}}
\newcommand{\Gr}{{\mbox {\rm Gr}}}
\newcommand{\Diff}{{\mbox {\rm Diff}}}
\newcommand{\Symp}{{\mbox {\rm Symp}}}
\newcommand{\symp}{{\mbox {\footnotesize{\rm symp}}}}
\newcommand{\Ct}{{\mbox {\rm Ct}}}
\newcommand{\Uns}{{\mbox {\rm Uns}}}
\newcommand{\rank}{{\mbox {\rm rank}}}
\newcommand{\sign}{{\mbox {\rm sign}}}
\newcommand{\supp}{{\mbox {\rm supp}}}
\newcommand{\Spin}{{\mbox {\rm Spin}}}
\newcommand{\Sp}{{\mbox {\rm sp-codim}}}
\newcommand{\Int}{{\mbox {\rm Int}}}
\newcommand{\Hom}{{\mbox {\rm Hom}}}
\newcommand{\codim}{{\mbox {\rm codim}}}
\newcommand{\ord}{{\mbox {\rm ord}}}
\newcommand{\Iso}{{\mbox {\rm Iso}}}
\newcommand{\corank}{{\mbox {\rm corank}}}
\def\mod{{\mbox {\rm mod}}}
\newcommand{\pt}{{\mbox {\rm pt}}}
\newcommand{\enP}{\hfill $\Box$ \par\vspace{5truemm}}
\newcommand{\qed}{\hfill $\Box$ \par}
\newcommand{\spe}{\vspace{0.4truecm}}
%\newcommand{\dfrac}[2]{\frac{\displaystyle{#1}}{\displaystyle{#2}}}
%%%%%%%%%%%%%%%%%

\title{Global classification of curves \\ 
on the symplectic plane} 

\author{Goo ISHIKAWA\thanks{\scriptsize 
Partially supported by 
Grants-in-Aid for Scientific Research, No. 14340020. }}

\renewcommand{\thefootnote}{\fnsymbol{footnote}}

\date{ }
\maketitle

\section{Introduction.}

In \cite{IJ}, 
we considered the local classification of 
plane curves on the symplectic plane. In particular, 
we introduced the number \lq\lq symplectic defect", 
which represents 
the difference of two natural 
equivalence relations on plane curves, the equivalence by diffeomorphisms 
and that by symplectomorphisms. 
For an immersion, two equivalence relations coincide, so 
the symplectic defect is null. For complicated singularities, 
the symplectic defects turn out to be positive. 

In this paper we consider the global symplectic classification problem.  
First we give the exact classification result under symplectomorphisms, for 
the case of generic plane curves, namely immersions 
with transverse self-intersections. 
Then, for a given diffeomorphism class of a generic plane curve, 
the set of symplectic classes form 
{\it the symplectic moduli space}
which we completely describe by its global topological term 
(Theorem \ref{generic immersion}). 
In the general plane curves with singularities, 
the difference between symplectomorphism and diffeomorphism classifications 
is clearly described 
by {\it local symplectic moduli spaces} of singularities 
and a global topological term. 
Thus, up to the classification by diffeomorphisms, 
the global problem is reduced to the local classification problem 
(Theorem \ref{of finite type}). 
We introduce the symplectic moduli space of a global plane curve and 
the local symplectic moduli space of a plane curve singularity  
as quotients of mapping spaces, and 
we endow them with differentiable 
structures in a natural way. 
Actually we treat {\it labelled} plane curves and 
labelled symplectic moduli spaces. For a plane curve, 
we label all compact domains 
surrounded by it and 
all singular points, and consider the classification problem 
of plane curves isotopic to the given plane curve by symplectomorphisms 
preserving the labelling.

\

Let $f : S^1 \to \R^2$ be a generic immersion of the circle $S^1$ 
in the symplectic plane $\R^2$ with the standard symplectic (area) form 
$\omega_0 = dx \wedge dy$. 
Clearly the areas of domains surrounded by the curve $f(S^1)$ 
are invariant under symplectomorphisms. Thus, denoting 
the first Betti number of $f(S^1)$ by $r$, we see 
the curves isotopic to $f$ 
have $r$-dimensional symplectic moduli. 

We denote by $C^\infty(S^1, \R^2)$ the 
space of $C^\infty$ mappings from $S^1$ to $\R^2$, 
which has the natural action (from \lq\lq right") 
of the group $\Diff^+(S^1)$ 
consisting of orientation-preserving diffeomorphisms on $S^1$. 
Thus $C^\infty(S^1, \R^2)/\Diff^+(S^1)$ 
denotes the space of {\it oriented} curves. 
The space $C^\infty(S^1, \R^2)/\Diff^+(S^1)$  has 
the action (from \lq\lq left") of the group $\Diff^+(\R^2)$ 
(resp. $\Symp(\R^2)$) consisting of 
orientation-preserving diffeomorphisms 
(resp. symplectomorphisms) on $\R^2$. 
For each oriented curve $f \in C^\infty(S^1, \R^2)/\Diff^+(S^1)$, 
we denote by $\Diff^+(\R^2)f$ the orbit through $f$ via the 
action of $\Diff^+(\R^2)$. Thus $\Diff^+(\R^2)f$ consists of 
oriented curves of form $\rho\circ f$ for orientation preserving 
diffeomorphisms $\rho$. 
Similarly the space $\Diff^+(\R^2)f/\Symp(\R^2)$ means the 
quotient space by the $\Symp(\R^2)$-action 
of $\Diff^+(\R^2)f$ in 
$C^\infty(S^1, \R^2)/\Diff^+(S^1)$. 
(The expression $\Symp(\R^2)\setminus\Diff^+(\R^2)f$ 
would be better, but I do not adopt it.) 

We call the quotient space $\Diff^+(\R^2)f/\Symp(\R^2)$ 
{\it the symplectic moduli space} of $f$ and denote it by 
${\mathcal M}_\symp(f)$. 
It describes the symplectic 
classification of a fixed isotopy class of an oriented plane curve.

To study the moduli space minutely, 
we label the $r$-domains surrounded by the curve 
$f(S^1)$ as $D_1, D_2, \dots, D_r$ for a 
$f \in C^\infty(S^1, \R^2)/\Diff^+(S^1)$. 
Then, for each $\rho \in \Diff^+(\R^2)$, we label bounded $r$-domains 
surrounded by $(\rho\circ f)(S^1)$ 
as $\rho(D_1), \rho(D_2), \dots, \rho(D_r)$ 
induced by the labelling for $f$. 
\begin{figure}[htbp]
  \begin{center}
\includegraphics[width=4truecm, height=4truecm, clip]{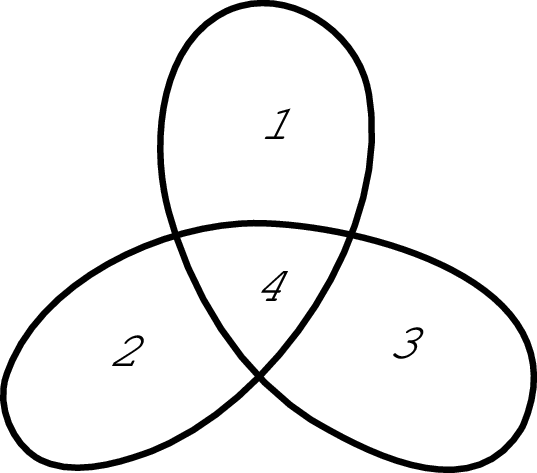} 
 %\vspace{-0.4truecm}
 \caption{A labelling of a generic plane curve.}
    \label{label3}
  \end{center}
 \end{figure}%
We set  
$$
\widetilde{\mathcal M}_\symp(f) = \Diff^+(\R^2)/\!\!\sim_f, 
$$
where 
we call $\rho, \rho'  \in \Diff^+(\R^2)$ are equivalent via $f$, and write 
 $\rho \sim_f \rho'$, 
if the exists a symplectomorphism $\tau$ such that 
 $\tau\circ\rho\circ f = \rho'\circ f$ up to $\Diff^+(S^1)$ and $\tau$ 
preserves
the given labelling: $\tau(\rho(D_j)) = \rho'(D_j), 1\leq j \leq r$. 
We call $\widetilde{\mathcal M}_\symp(f)$ the {\it 
labelled symplectic moduli space} of $f$. 
Note that $\widetilde{\mathcal M}_\symp(f)$ does not 
depend on the chosen labelling of $f$. 

The natural projection $\pi : 
\widetilde{\mathcal M}_\symp(f) \to {\mathcal M}_\symp(f)$ 
is defined by $\pi([\rho]) = [\rho\circ f]$ and $\pi$ 
is a finite (at most $r!$, $r$ factorial) to one. 

Then in this paper we show in fact:  

\bet
\label{generic immersion}
If $f \in C^\infty(S^1, \R^2)/\Diff^+(S^1)$ 
is a generic 
%labelled 
immersion, 
then the labelled symplectic moduli space $\widetilde{\M}_\symp(f)$ 
is diffeomorphic to the relative cohomology space  
$H^2(\R^2, f(S^1), \R) \cong \R^r$. 
%(\cong H^1(f(S^1), \R))$
\ent
%$$
%\Diff^+(\R^2)f/\Symp(\R^2) \quad
%{\mbox{\Large$\cong$}}_{\mbox{\rm{\small{diffeo}}}} 
%\quad H^1(f(S^1), \R). 
%$$ 

The labelled symplectic moduli space 
$\widetilde{\mathcal M}_\symp(f)$ 
has a canonical differentiable structure. 
We claim in Theorem \ref{generic immersion} that 
the labelled symplectic moduli space of $f$ 
with the differentiable structure is {\it diffeomorphic} to 
$\R^r$, $r = \dim_{\R}H^2(\R^2, f(S^1), \R) = 
\dim_{\R}H_1(f(S^1), \R)$.

Actually we are going to give a diffeomorphism 
between 
%$\Diff^+(\R^2)f/\Symp(\R^2)$ 
$\widetilde{\mathcal M}_\symp(f)$ 
and the positive cone in $H^2(\R^2, f(S^1), \R)$. 
Note that the relative cohomology group $H^2(\R^2, f(S^1), \R)$ 
over $\R$ 
is isomorphic to the vector space $H_2(\R^2, f(S^1), \R)^* 
= \Hom_{\R}(H_2(\R^2, f(S^1), \R), \R)$.  
The orientation of $\R^2$ and labelling of the bounded domains 
surrounded by $f(S^1)$ give the canonical basis $[D_1], [D_2], \dots, 
[D_r]$ of 
$H_2(\R^2, f(S^1), \R)$. 
The positive cone $H^2(\R^2, f(S^1), \R)_{>0}$
%of $H^2(\R^2, f(S^1), \R) \cong H_2(\R^2, f(S^1), \R)^*$ 
is defined 
by 
$$
H^2(\R^2, f(S^1), \R)_{>0} = 
\{ \alpha \in H^2(\R^2, f(S^1), \R) \ \mid \ \alpha([D_j]) > 0, \ 
1 \leq j \leq r 
\}. 
$$

The diffeomorphism of $\widetilde{\mathcal M}_\symp(f)$
%$\Diff^+(\R^2)f/\Symp(\R^2)$ 
and 
$H^2(\R^2, f(S^1), \R)_{>0}$ is given actually by the mapping 
$$
\varphi : 
\widetilde{\mathcal M}_\symp(f)
%\Diff^+(\R^2)f/\Symp(\R^2) 
\to H^2(\R^2, f(S^1), \R)_{>0}, 
$$
defined by 
$$
\varphi : [\rho] \mapsto \left( [D_j] 
%\in H_2(\R^2, f(S^1), \R) 
\mapsto \int_{\rho(D_j)}\omega_0 = \int_{D_j}\rho^*\omega_0 
\right). 
$$
$1 \leq j \leq r, \ \omega_0 = dx\wedge dy$. 

\smallskip 

The symplectic moduli space ${\mathcal M}_\symp(f) = 
\Diff^+(\R^2)f/\Symp(\R^2)$ 
is obtained as a quotient of $\widetilde{\mathcal M}_\symp(f)$. 
A {\it symmetry} of a generic immersion $f : S^1 \to \R^2$ 
is an orientation preserving diffeomorphism $\rho : \R^2 \to \R^2$ 
such that $\rho\circ f = f\circ \sigma$ for some $\sigma \in 
\Diff^+(S^1)$. 
We denote by $S_f$ the group of symmetries of $f$. 
Then $S_f$ induces
a subgroup $G_f$ of the permutation group $S_r$ of the 
$r$-bounded domains of $\R^2 \setminus f(S^1)$. 
Then $G_f$ naturally acts on 
$H_2(\R^2, f(S^1), \R)$ and on $H_2(\R^2, f(S^1), \R)^* \cong 
H^2(\R^2, f(S^1), \R) \cong \R^r$ as permutation of coordinates. 

\bee
Let $f$ be a generic immersion as in Figure \ref{label3}: 
Then $G_f \subset S_4$ is a cyclic group of order $3$.  
\ene

By Theorem \ref{generic immersion}, we have the following: 

\bec
\label{generic immersion2}
The symplectic moduli space ${\mathcal M}_\symp(f)$ is diffeomorphic to 
the finite quotient $\R^r/G_f$ of $\R^r$. 
\enc

In fact, the action of $G_f \subset S_r$ commutes with the diffeomorphism 
$(\R_{>0})^r \to \R^r$ defined by $(x_1, \dots, x_r) 
\mapsto (\log x_1, \dots, \log x_r)$.

\

We denote by ${\mathcal M}_\symp$ the whole orbit space of 
$C^\infty(S^1, \R^2)$ by the $\Diff^+(S^1)\times \Symp(\R^2)$-action 
(the right-left-symplectic action): 
$$
{\mathcal M}_\symp := 
C^\infty(S^1, \R^2)/(\Diff^+(S^1)\times \Symp(\R^2)). 
$$
Note that ${\mathcal M}_\symp$ is a non-Hausdorff space, with respect to 
the quotient topology of $C^\infty$ topology. 
The non-Hausdorffness comes from the adjacencies of 
$\Diff^+(\R^2)\times\Diff^+(S^1)$-orbits. 

Besides the singular topology, we have the decomposition 
$$
{\mathcal M}_\symp = \bigcup {\mathcal M}_\symp(f), 
$$
where $f$ runs over representatives  
of the set of isotopy types 
%$C^\infty(S^1, \R^2)/(\Diff^+(\R^2)\times \Diff^+(S^1))$ 
of oriented plane curves. 
We ask the structure of 
each stratum ${\mathcal M}_\symp(f)$ itself. 
Then Corollary \ref{generic immersion2} guarantees that 
each \lq\lq open stratum" ${\mathcal M}_\symp(f)$ 
is a finite quotient of an affine space 
where $f$ is a generic immersion. 

\smallskip

The method to provide a \lq\lq differentiable structure\rq\rq 
to a mapping space quotient (a moduli space) 
should be not unique \cite{Eells}\cite{Omori}. 
For instance, consider the problem how to define 
a differentiable structure on a mapping space $C^\infty(N, M)$ itself 
for $C^\infty$ manifolds $N$ and $M$.  
Then one of the standard methods 
seems to define, first, Fr{\' e}chet differentiable functions on 
the Banach manifolds $C^r(N, M)$, for each finite $r$, and 
regard $C^\infty(N, M)$ as the inverse limit of $C^r(N, M)$ 
to define the structure sheaf of differentiable functions on it. 
However in this paper we apply another method: 
We regard Fr{\' e}chet differential of a functional 
as a kind of \lq\lq total differential\rq\rq. 
Then we could consider, instead, \lq\lq partial differentials\rq\rq. 
Namely, to define differentiable functions on $C^\infty(N, M)$, 
first we define the notion of finite dimensional differentiable families 
in $C^\infty(N, M)$ by the very classical and natural manner. 
Then we call a function on $C^\infty(N, M)$ differentiable 
if its restriction to any finite dimensional 
family in $C^\infty(N, M)$ is of class $C^\infty$ in the ordinary sense. 
See \S \ref{Diff str.}. 

\

Theorem \ref{generic immersion} 
is generalised to more singular curves. 
To state the generalisation, first we treat 
the local case. 

A multi-germ 
$f = f_{y_0} : (S^1, S) \to (\R^2, y_0)$ 
at a finite set $S \subset S^1$ 
is called {\it of finite codimension} (or ${\mathcal A}$-{\it finite}) 
in the sense of Mather 
if $f_{y_0}$ is determined by its finite jet up to diffeomorphisms 
(or ${\mathcal A}$ equivalence). See
\cite{Mather3}\cite{Wall}. 
We denote by $\Diff^+_0(S^1, S)$ the group of 
orientation preserving diffeomorphism-germs $(S^1, S)$ 
fixing $S$ pointwise, and 

we treat $f_{y_0}$ up to 
$\Diff_0^+(S^1, S)$. 
Then the local image of $f_{y_0}$ divides $(\R^2, 0)$ into 
several domains. We label them. 
Then, for any orientation-preserving diffeomorphism-germ $\rho 
\in \Diff^+(\R^2, y_0)$, the labelling of $\rho\circ f$ 
is induced. 
Two diffeomorphism-germs $\rho, \rho' \in 
\Diff^+(\R^2, y_0)$ are equivalent via $f_{y_0}$, and write 
$\rho \sim_f \rho'$, if 
there exists a symplectomorphism-germ $\tau \in 
\Symp(\R^2, 0)$ 
such that $\tau\circ\rho\circ f_{y_0} = \rho'\circ f_{y_0}$ up to 
$\Diff_0^+(S^1, S)$ 
and $\tau$ preserves the labelling. 
Thus we define the {\it local labelled symplectic moduli space} 
by 
$$
\widetilde{\mathcal M}_\symp(f_{y_0}) := 
\Diff^+(\R^2, y_0)/\sim_f. 
$$
Moreover we define the {\it local symplectic moduli space} by 
$$
{\mathcal M}_\symp(f_{y_0}) := 
\Diff^+(\R^2, y_0)f_{y_0}/\Symp(\R^2, y_0). 
$$

Note that the space of map-germs 
$$
C^\infty((N, S), (M, y_0)) := 
\{ f : (N, S) \to (M, y_0) \ C^\infty {\mbox{\rm \ map-germs}} \} 
$$
is a quotient space of $C^\infty(N, M)$, so also 
it has the differentiable structure. 

In particular $\widetilde{\mathcal M}_\symp(f_{y_0})$ and 
${\mathcal M}_\symp(f_{y_0})$ are mapping space quotient 
have natural differentiable structures. Moreover 
there exists the canonical projection 
$\pi : \widetilde{\mathcal M}_\symp(f_{y_0}) 
\to {\mathcal M}_\symp(f_{y_0})$ defined by 
$\pi([\rho]) = [\rho\circ f]$ modulo 
$\Diff_0^+(S^1, S)$. 

\ 

Now returning to the global case, we consider again an oriented curve 
$f : S^1 \to \R^2$ up to $\Diff^+(S^1)$, namely, 
$f \in C^\infty(S^1, \R^2)/\Diff^+(S^1)$. 
Then we call $f$ {\it of finite type} if,  
for some 
(and for any) representative $f : S^1 \to \R^2$ of 
$f \in C^\infty(S^1, \R^2)/\Diff^+(S^1)$, 
except for a finite number of 
points $y_0 \in f(S^1)$, 
the multi-germ 
$f_{y_0} : (S^1, f^{-1}(y_0)) \to (\R^2, y_0)$ is a stable multi-germ, 
namely a single immersion-germ or a transversal two-immersion-germ, 
and, even if  $f_{y_0}$ is unstable, 
$f^{-1}(y_0)$ is a finite set in $S^1$ and $f_{y_0}$ 
is of finite codimension. 
The condition means roughly that the $\Diff^+(\R^2)$-orbit 
through $f$ in $C^\infty(S^1, \R^2)/\Diff^+(S^1)$ 
is of finite codimension.

\begin{figure}[htbp]
  \begin{center}
\vspace{0.5truecm}
 \includegraphics[width=4truecm, height=4truecm, clip]{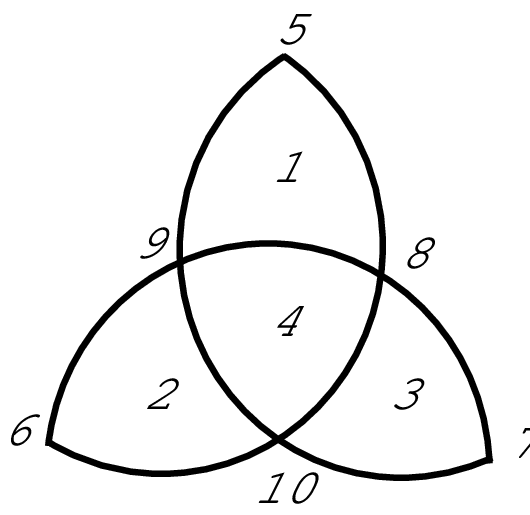} 
 %\vspace{-0.8truecm}
 \caption{A labelling of a plane curve of finite type.}
    \label{label4}
  \end{center}
 \end{figure}

If $f \in C^\infty(S^1, \R^2)/\Diff^+(S^1)$ is of finite type, 
then $f(S^1)$ divides $\R^2$ into a finite number of bounded 
domains and one unbounded domain. 
Then we define the the {\it labelling} of $f$ 
as the labelling of bounded domains $D_1, \dots, D_r$ 
and the multiple or singular values $y_1, \dots, y_s$ of 
$f$ in $\R^2$ (Figure \ref{label4}, where $r = 4, s = 6$). 
We define, 
similarly to the case of generic immersions, 
the {\it labelled symplectic moduli space} 
of plane curve $f$ of 
finite type by 
$$
\widetilde{\mathcal M}_\symp(f) := 
\Diff^+(\R^2)/\sim_f,
$$
where 
$\rho \sim_f \rho'$ if 
$\tau\circ\rho\circ f = \rho'\circ f$ for some 
$\tau \in \Symp(\R^2)$ preserving the labellings induced by $\rho$ and 
$\rho'$. 

Moreover we define 
the {\it symplectic moduli space} of plane curve $f$ of 
finite type by 
$$
{\mathcal M}_\symp(f) := \Diff^+(\R^2)f/\Symp(\R^2). 
$$

\bet
\label{of finite type}
{\rm (Localisation Theorem)} 
If $f \in C^\infty(S^1, \R^2)/\Diff^+(S^1)$ 
is of finite type, then 
$$
\widetilde{\mathcal M}_\symp(f) 
\ \cong_{\rm{diffeo.}} 
\prod_{y_0 \in f(S^1)} 
\widetilde{\mathcal M}_\symp(f_{y_0}) 
\times \R^r, 
$$
where $r = \dim_\R H^2(\R^2, f(S^1), \R)$. 
Moreover 
$$
{\mathcal M}_\symp(f) 
\ \cong_{\rm{diffeo.}} 
\left(\prod_{y_0 \in f(S^1)} 
\widetilde{\mathcal M}_\symp(f_{y_0}) 
\times \R^r \right)/G_f, 
$$
where $G_f \subset S_{r'}$ is the group induced by 
the symmetry group of $f$, $r'$ being 
$r$ plus the number of unstable singular values of $f$. 
\ent 

Note that  
$\widetilde{\mathcal M}_\symp(f_{y_0})$ is just a point 
if $f_{y_0}$ is stable, namely if it is a single immersion-germ 
or a transverse self-intersection. Therefore the product in Theorem 
\ref{of finite type} turns out to be a finite product. 

Theorem \ref{of finite type} can be regarded as the \lq\lq localisation 
theorem\rq\rq for the global labelled moduli space of the isotopy type of 
a singular plane curve. 

\vspace{0.2truecm}

The diffeomorphism between $\widetilde{\mathcal M}_\symp(f)$  and 
the product of local symplectic moduli spaces and an open cone of 
$H^2(\R^2, f(S^1), \R)_{>0}$ is given by the mapping 
$$
\Phi : \widetilde{\mathcal M}_\symp(f) \to 
\left(\prod_{y_0 \in f(S^1)} 
\widetilde{\mathcal M}_\symp(f_{y_0}) \right)\times 
H^2(\R^2, f(S^1), \R)_{>0} 
$$
defined by 
$$
\Phi([\rho]) := 
\left( ([\eta_{\rho(y_0)}\circ 
\rho])_{y_0 \in f(S^1)}, 
\  \varphi([\rho]) \right), 
$$
where $\eta_{\rho(y_0)} : (\R^2, \rho(y_0)) \to (\R^2, y_0)$ 
is {\it any} symplectomorphism-germ. Note that 
$[\eta_{\rho(y_0)}\circ \rho] \in \widetilde{\mathcal M}(f_{y_0})$ 
does not depend on the choice of $\eta_{\rho(y_0)}$.

\ 

Here are some examples of local (labelled) symplectic moduli spaces: 

\bee
{\rm 
(1) For the germ $f : (\R, 0) \to (\R^2, 0)$ of type $A_2$ 
defined by $f(t) = (t^2, t^3)$ 
we have $\widetilde{\mathcal M}_\symp(f) = {\mathcal M}_\symp(f) = 
\{ \pt\}$. 

(2) For multi-germ $f : (\R, \{ 0, 1\}) \to 
(\R^2, 0)$ (resp. $f : (\R, \{ 0, 1\}) \to 
(\R^2, 0)$; $f : (\R, \{ 0, 1, 2\}) \to 
(\R^2, 0)$) defined by 
$f(t) = (t, 0)$ near $t = 0$ and $f(s) = (0, s-1)$ near $s = 1$ 
(resp. 
$f(t) = (t, 0)$ near $t = 0$ and $f(s) = (0, (s-1)^2)$ near $s = 1$; 
$f(t) = (t, 0)$ near $t = 0$, $f(s) = (0, s-1)$ near $s = 1$ and 
$f(r) = (r-2, r-2)$ near $r = 2$)), we have 
$\widetilde{\mathcal M}_\symp(f) = {\mathcal M}_\symp(f) = \{ \pt\}$. 

(3)
Let $f : (\R, 0) \to (\R^2, 0)$ be a map-germ 
of type $E_{12}$ defined by $f(t) = 
(t^3, t^7)$. Then $f$ has the symplectic normal form 
$(t^3, t^7 + \lambda t^8)$. 
and we see $\widetilde{\mathcal M}_\symp(f) 
= {\mathcal M}_\symp(f)$ is diffeomorphic to $\R$. 

(4)
Let $f : (\R, 0) \to (\R^2, 0)$ be a map-germ 
of type $W_{18}$ defined by $f(t) = 
(t^4, t^7 + t^9)$. 
Then $f$ has the symplectic normal form 
$(t^4, t^7 + \lambda_1t^9 + \lambda_2t^{13}), \lambda_1 > 0$, 
and we see 
$\widetilde{\mathcal M}_\symp(f) 
= {\mathcal M}_\symp(f)$ is diffeomorphic to $\R^2$. 

(5) 
Let $f : (\R, 0) \to (\R^2, 0)$ be a map-germ 
of type $E_{24}$ defined by $f(t) = 
(t^3, t^{13} + t^{14})$. 
Then $f$ has the symplectic normal form 
$(t^3, t^{13} + \lambda_1t^{14} + \lambda_ 2t^{17} + \lambda_3 t^{20}), \lambda_1 > 0$, 
and we see 
$\widetilde{\mathcal M}_\symp(f) 
= {\mathcal M}_\symp(f)$ is diffeomorphic to $\R^3$. 
}
\ene

As an application of Theorem \ref{of finite type}, we have for 
instance: 

\bee
{\rm 
(cf. \cite{IJ}). 
If $f : S^1 \to \R^2$ is an injective $C^\infty$ map 
with just one singular point, where $f$ 
is locally diffeomorphic to $t \to (t^3, t^{13} + t^{14})$ (of type $E_{24}$). 
Then the local symplectic 
moduli space for the singular point is 
diffeomorphic to $\R^3$ and the symplectic moduli space 
${\mathcal M}_\symp(f) = \Diff^+(\R^2)f/\Symp(\R^2)$ of $f$ 
is diffeomorphic to 
$\R^3 \times \R \cong \R^4$. 
}
\ene 

\ber
{\rm 
In all known examples, we observe that the local symplectic moduli 
space are differentiable manifolds and 
the symplectic defect is interpreted as the dimension of the tangent space to the moduli space. 
We conjecture that this holds in general for multi-germs $f$ of finite type. 
Moreover,  
we conjecture that $\widetilde{\mathcal M}_\symp(f)$ is 
diffeomorphic to $\R^\ell$ with 
$\ell = 
\dfrac{1}{2}\mu(f) - \codim(f), 
$
where $\mu(f)$ is Milnor number of $f$ and 
$\codim(f)$ is the 
${\mathcal A}_e$-codimension of $f$ (\cite{Wall}). 

Over the complex numbers, we conjecture that similarly defined 
${\mathcal M}_\symp(f)$ is a $K(\pi, 1)$ space and 
the universal covering of ${\mathcal M}_\symp(f)$ 
is diffeomorphic to $\C^\ell$ with $\ell = 
\dfrac{1}{2}\mu(f) - \codim(f)$. See \cite{Varchenko}. 
}
\enr

\ 

Similar result to Theorem \ref{of finite type} holds for curves with finite components: 
Let $N$ be a one dimensional closed manifold, i.e. a finite disjoint union of 
circles. 

\bet
\label{multi-components}
If $f \in C^\infty(N, \R^2)/\Diff^+(N)$ 
is of finite type, then, for the labelled symplectic moduli space, 
we have 
$$
%\dfrac{\Diff^+(\R^2)f}{\Symp(\R^2)} 
\widetilde{\mathcal M}_\symp(f) 
\quad
{\mbox{\large$\cong$}}_{\mbox{\rm{\small{diffeo.}}}} \  
\left(
{\prod}_{y_0 \in f(N)} \widetilde{\mathcal M}_\symp(f_{y_0}) 
\right) 
\times H^2(\R^2, f(N), \R). 
$$ 
Moreover ${\mathcal M}_\symp(f)$ 
is diffeomorphic to $\widetilde{\mathcal M}_\symp(f)/G_f$. 
\ent

Furthermore, Theorems \ref{of finite type} and \ref{multi-components} 
are generalised to 
curves with a finite number of components in a symplectic surface: 

A two dimensional symplectic manifold is called 
a {\it symplectic surface}. A symplectic surface $(M, \omega)$ {\it has 
a bounded area} (resp. {\it an unbounded area}) if 
$\displaystyle{\int_M\omega < \infty}$ (resp. $ = \infty$). 
Denote by $\Diff_0(M)$ (resp. $\Symp_0(M)$) 
the group of diffeomorphisms on $M$ isotopic to the identity 
(resp. symplectomorphisms on $M$ isotopic to the identity 
through symplectomorphisms.) 
Then we set 
$$
\widetilde{\mathcal M}_\symp(f) := 
\Diff_0(M)/\!\!\sim_f, 
$$
where $\rho \sim_f \rho'$ if there exists $\tau \in 
\Symp_0(f)$ preserving a labelling of $f$ and 
satisfying that $\tau\circ\rho\circ f = \rho'\circ f$. 
Moreover we set 
$$
{\mathcal M}_\symp(f) := \Diff_0(M)f/\Symp_0(f). 
$$
Then we have 

\bet
\label{surface}
Let $(M, \omega)$ be a connected symplectic surface, 
$N$ a one dimensional closed manifold and 
$f \in C^\infty(N, M)/\Diff^+(N)$ 
an oriented curve of finite type. Denote by $r$ the 
number of connected components of $M \setminus f(N)$ 
with bounded areas. 

{\rm (1)} 
If $(M, \omega)$ itself has an unbounded area, then 
$$
\widetilde{\mathcal M}_\symp(f) \quad
{\mbox{\large$\cong$}}_{\mbox{\rm{\small{diffeo.}}}}  \ 
\left(
{\prod}_{y_0 \in f(S^1)} \widetilde{\mathcal M}_\symp(f_{y_0}) 
\right) 
\times \R^{r}. 
$$ 

{\rm (2)} 
If $(M, \omega)$ itself has a bounded area, then 
$$
\widetilde{\mathcal M}_\symp(f) \quad
{\mbox{\large$\cong$}}_{\mbox{\rm{\small{diffeo.}}}} \  
\left(
{\prod}_{y_0 \in f(S^1)} \widetilde{\mathcal M}_\symp(f_{y_0}) 
\right) 
\times \R^{r-1}. 
$$ 

Moreover, in any case, ${\mathcal M}_\symp(f)$ 
is diffeomorphic to $\widetilde{\mathcal M}_\symp(f)/G_f$. 
\ent

\ber
{\rm 

(1) If $(M, \omega)$ has a bounded area, then the 
sum of areas of domains surrounded by the curve 
must be equal to the total area. This restriction reduces 
the dimension of the moduli space by one. 

(2) If $M$ is a closed surface, then 
$(M, \omega)$ has a bounded area for any symplectic form $\omega$. 

(3) If $(M, \omega) = (\R^2, \omega_0)$, the standard symplectic $\R^2$, 
then 
$r$ is equal to  $\dim_{\R}H^2(\R^2, f(N), \R)$. 
}
\enr

\

To describe completely the symplectic classification, 
we need to solve also the classification problem of plane curves 
by diffeomorphisms. See \cite{Arnold} for the global classification 
of generic immersions by diffeomorphisms (or homeomorphisms). 
For the local diffeomorphism classification (over the complex), 
there are several detailed studies 
\cite{Ebey}\cite{Zariski}\cite{Wahl}\cite{BG}\cite{Teissier}\cite{LaP}. 
Note that \cite{BG} gives the classification of simple singularities of 
parametric plane curves. 
The classification results are improved by \cite{LuP}\cite{HH}\cite{IJ3} 
recently. 

Note that, in \cite{IJ3}, we give the 
diffeomorphism classification of unimodal singularities 
of parametric plane curves, and moreover 
we describe concretely the symplectic 
moduli space of simple and unimodal plane curve singularities 
in the complex analytic category.

On the plane, the classification by symplectomorphisms 
is of course same as the classification by area-preserving diffeomorphisms. 
Nevertheless our classification problem seems to belong to symplectic geometry 
not to volume geometry, 
from both the motivation and the mathematical reason. 
In fact, the diffeomorphism invariance of symplectic defects for plane curves
is naturally generalised to higher dimensional cases 
in term of Lagrangian varieties in symplectic spaces(\cite{IJ2}). 
See several related results \cite{Varchenko}\cite{Lando}\cite{Garay}\cite{DR}.

Naturally we can ask the similar results to Theorems 
of this paper, for higher dimensional cases: 
The symplectic classification of Lagrangian surfaces 
of a given Lagrangian isotopy class in $\R^4$ for instance. 
For that, we have to clarify on the topology (or symmetries) 
of generic Lagrangian immersions 
as well as the local classification of singularities of Lagrangian 
surfaces (\cite{Givental}\cite{Inv}). 
Both problems seems to be very interesting.

Also, setting various geometric structures 
on the symplectic moduli spaces 
like in \cite{Hitchin1}\cite{Hitchin2}\cite{Brylinski} 
would be very interesting problem, which is still open as far as I know.

\

In \S \ref{Diff str.}, we define a differentiable structure an any 
mapping space quotient. 

In \S \ref{Mild spaces.}, we introduce a class of mapping space quotients 
which is reasonable to study for our purpose. 

In 
\S \ref{Differentiable structure on a moduli space.}, 
we define differentiable structures on the (labelled) symplectic moduli space. 
The symplectic classification problem of curves of fixed diffeomorphism 
class can be translated to the classification of forms by diffeomorphisms 
fixing a curve (\cite{Zhitomirskii}\cite{Zhitomirskii2}\cite{DJZ}\cite{DR}). 
In \S \ref{The space of symplectic forms.}, 
we relate the space of symplectic forms with 
the symplectic moduli spaces. 
Then we give proofs of main theorems. 

\ 

The starting point of the present paper is 
the symplectic bifurcation problem studied by the joint work \cite{IJ} 
with S. Janeczko. I would like to thank him for his continuous 
encouragement.  
I am grateful to H. Sato and S. Izumi for their valuable comment.

\section{Differentiable structure on a mapping space quotient.} 
\label{Diff str.} 

We denote by $C^\infty(N, M)$ the space of $C^\infty$ 
mappings from a (finite dimensional) $C^\infty$ manifold $N$ to 
a (finite dimensional) $C^\infty$ manifold $M$. 
In this section, also $P, Q, L, K$ always 
designate (finite dimensional) $C^\infty$ manifolds respectively.

Let 
$X \subseteq C^\infty(N, M)$ be a subset. Then, such a set 

$X$ is a {\it mapping space}. 
Let $X/\!\!\sim$ be any quotient of $X$ 
under an equivalence relation $\sim$ on $X$. 
We give on the quotient space $X/\!\!\sim$ the quotient topology 
of ${X}$ with the relative topology 
of the $C^\infty$ topology on $C^\infty(N, M)$, 
not Whitney (fine) $C^\infty$ topology.  
Such space is called a {\it mapping space quotient}. 
Then we will endow, in the following five steps,  
a differentiable structure 
with the mapping space quotient $X$, 
depending just on the representation 
$X/\!\!\sim \ \leftarrow X \subseteq C^\infty(N, M)$. 
We note that the notion of differentiable structures  
can provided just by defining the notion of diffeomorphisms. 
Therefore our goal is to define the notion of \lq\lq diffeomorphisms\rq\rq. 

(i) 
We call a mapping 

$h : P \to X$  {\it differentiable} (or $C^\infty$) if 
there exists a $C^\infty$ mapping (between manifolds) 
$H : P\times N \to M$ 
satisfying $H(p, x) = h(p)(x) \in M$,  
($p \in P, x \in N$). 

(ii) 

We call a mapping 

$k : X \to Q$  {\it differentiable} if 

$k$ is a continuous mapping and, 

for any differentiable mapping 

$h : P \to X$ in the sense of (i), the composition 

$k\circ h : P \to Q$ 
is a $C^\infty$ mapping between manifolds.

Now, if 
$\sim$ is an equivalence relation on a mapping space 

$X$, then we get 

 the quotient space $X/\!\!\sim$. 
Then the canonical projection 

$\pi = \pi_X : X \to X/\!\!\sim$ is defined by $\pi(x) = [x]$ (the equivalence class 

of $x$). 

(iii) 

We call a mapping 
$\ell : X/\!\!\sim \ \to Q$ {\it differentiable} if 
the composition $\ell\circ\pi : X \to Q$ with the projection $\pi$ 

is differentiable in the sense of (ii).

(iv)

We call 
a mapping 
$\varphi : X/\!\!\sim \ \to Y/\!\!\approx$ 

from a mapping space quotient $X/\!\!\sim$ 
to another mapping space quotient 

$Y/\!\!\approx \ \leftarrow 
Y \subseteq C^\infty(L, K)$ {\it differentiable} if 
$\varphi$ is a continuous mapping 

and, for any open subset $U \subseteq Y/\!\!\approx 
(\leftarrow \pi_Y^{-1}(U) \subseteq C^\infty(L, K))$ and 
for any differentiable mapping 

$\ell : U \ 
\to Q$ in the sense of (iii), the composition 

$\ell\circ\varphi : \varphi^{-1}(U) (\leftarrow \pi_X^{-1}(\varphi^{-1}(U)) 
\subseteq C^\infty(N, M)) 
\to Q$ is 
differentiable in the sense of (iii).

(v) 

We call a mapping 
$\varphi : X/\!\!\sim \ \to Y/\!\!\approx$ a {\it diffeomorphism} 
if $\varphi$ is differentiable in the sense of (iv)C$\varphi$ is a bijection and 
the inverse mapping $\varphi^{-1} : Y/\!\!\approx \ \to X/\!\!\sim$ is also differentiable in the sense of (iv). 

Moreover 

we call two mapping space quotients 
$X/\!\!\sim$ and $Y/\!\!\approx$ {\it diffeomorphic} 

if there exists a diffeomorphism 
$\varphi : X/\!\!\sim \ \to Y/\!\!\approx$ in the sense of (iv). 

\

Now we give several related results: 
First, form the definition above, 
we immediately have that the differentiability is a local notion: 

\bel
A mapping $\varphi : X/\!\!\sim \ \to Y/\!\!\approx$ 
is a differentiable mapping 
if and only if $\varphi$ is locally a differentiable mapping, 
namely, if and only if, for any $x_0 \in X/\!\!\sim$, 
there exists a neighborhood 
$U$ of $x_0$ in $X/\!\!\sim$, such that $\varphi\vert_U : U 
\to Y/\!\!\approx$ is a differentiable mapping. 
\enl

Also we observe the followings: 

\bel
For any $C^\infty$ manifold $P$, 
$P$ is diffeomorphic to $C^\infty(\{ \pt\}, P)$. 
\enl

\Proof
In fact the mapping $\varphi : P \to C^\infty(\{ \pt\}, P)$ 
defined by $\varphi(p)(\pt) = p, (p \in P)$ is a diffeomorphism. 
\QED

\bee
{\rm
Let $\QQ \subset \R = C^\infty(\pt, \R)$ be 
the set of rational numbers. 
Then a mapping $h : P \to \QQ$ from a $C^\infty$ manifold 
is differentiable if and only if $h$ is continuous (i.e. locally constant). 
A mapping $k : \QQ \to Q$ to a $C^\infty$ manifold is 
differentiable if and only if $k$ is continuous. 
}
\ene

We show useful lemmata which follow from the definition:

\bel
\label{Lemma1}
If $h : P \to X$ is differentiable in the sense of {\rm (i)}, then 
$\pi\circ h : P \to X/\!\!\sim$ is differentiable in the sense of 
{\rm (iv)}. 

\enl

\Proof
For any differentiable mapping $\ell : U (\subseteq X/\!\!\sim)  \to Q$ 

in the sense of (iii), 

the composition 

$\ell\circ \pi : \pi_X^{-1}(U) \to Q$ is differentiable in the sense of (ii). 

Therefore 

$(\ell\circ \pi)\circ h = 
\ell\circ(\pi\circ h) : (\pi\circ h)^{-1}(U) \to Q$ 
is differentiable. Hence $\pi\circ h$ differentiable in the sense of (iv). 
\QED

\bel
\label{diff-equiv}
The following two conditions are equivalent to each other:

{\rm (1)} $\varphi : X/\!\!\sim \ \to Y/\!\!\approx$ is 
differentiable in the sense of {\rm (iv)}. 

{\rm (2)} $\varphi : X/\!\!\sim \ \to Y/\!\approx$ is a continuous mapping 
and, for any differentiable mapping $h : P \to X$ in the sense of {\rm (i)}, 
$\varphi\circ\pi\circ h : P \to Y/\!\!\approx$ is 
differentiable in the sense of {\rm (iv)}. 
\enl

\Proof
(1) $\Rightarrow$ (2): 

Let 
$\ell : U (\subseteq Y/\!\!\approx) \to Q$ be a differentiable mapping in the sense of (iii). 

By (1), $\ell\circ(\varphi\circ\pi\circ h) 
= (\ell\circ\varphi\circ\pi)\circ h 
: (\varphi\circ\pi\circ h)^{-1}(U) \to Q$ is 
differentiable. 
Therefore $\varphi\circ\pi\circ h : P \to Y/\!\!\approx$ is  
a differentiable mapping in the sense of (iv). 

(2) $\Rightarrow$ (1): 

Let  

$\ell : U (\subseteq Y/\!\!\approx) \to Q$ be 
differentiable in the sense of (iii), and 
$h : P \to (\varphi\circ\pi)^{-1}(U)$ differentiable in the sense of (i). 
Then $h : (\varphi\circ\pi\circ h)^{-1}(U) 
\to (\varphi\circ\pi)^{-1}(U)$ is differentiable in the sense of (i). 
By (1), we have that $\ell\circ\varphi : \varphi^{-1}(U) \to Q$ 
differentiable in the sense of (iii). 
Therefore we see 
$(\ell\circ\varphi\circ\pi)\circ h = \ell\circ(\varphi\circ\pi\circ h) 
: P \to Q$ is $C^\infty$.  Therefore 
$\ell\circ\varphi\circ\pi : (\varphi\circ\pi)^{-1}(U) \to Q$ 
is differentiable, so is $\ell\circ\varphi : \varphi^{-1}(U) \to Q$ 
in the sense of (iii). 
Hence $\varphi$ is differentiable in the sense of (iv). 
\QED

\ber
{\rm
We remark that, in (2) of Lemma \ref{diff-equiv}, 
we use differentiable functions $h : P \to X$ to 
$X$, not to the quotient $X/\!\!\sim$, as \lq\lq test mappings", 
to teat the differentiability of $\varphi$. 
Actually, the class of differentiable mappings $P \to X/\!\!\sim$ 
depends heavily on the nature of the equivalence relation $\sim$. 
}
\enr

\bel
\label{cont}
A differentiable mapping $h : P \to X \subseteq C^\infty(N, M)$ in the sense of 

{\rm (i)} is a continuous mapping. 
\enl

\Proof

By the assumption, there exists a differentiable 
mapping $H : P \times N \to M$ 
which satisfies $H(p, x) = h(p)(x)$.  
Take an open subset of $C^\infty(N, M)$ of the form $W(r, K, U)$, 

where $K \subseteq N$ is a compact subset and $U \subseteq J^r(N, M)$ is an open subset.

Suppose, for a 
$p_0 \in P$,  

$h(p_0) = H\vert_{{p_0}\times N} : 
N \times M$ belongs to $W(r, K. U)$. 

Define $j^r_1H : P\times N \to J^r(N, M)$ by 
$j^r_1H(p, x) = j^r(H\vert_{p\times N})(x)$. 

Then $j^r_1H$ is a differentiable mapping in the ordinary sense. 

In particular it is continuous. From the assumption 

$h(p_0) \in W(r, K, U)$,  
$(j^r_1H)^{-1}(W(r, K, U))$ is an open neighborhood of $p_0\times K$. 

Since $K$ is compact, there exists an open neighborhood $V$ of $p_0$ 

such that $V\times K \subseteq (j^r_1H)^{-1}(W(r, K, U))$. 
This means that $p_0 \in V \subseteq h^{-1}(W(r, K, U))$. 

Therefore $h^{-1}(W(r, K, U))$ is open. 
Noting that 
$h^{-1}(W(r, K, U) \cap W(r', K', U')) = 
h^{-1}(W(r, K, U)) \cap h^{-1}(W(r', K', U'))$ and that 
$h^{-1}(\cup W_\nu) = \cup h^{-1}(W_\nu)$, we see that 
$h$ is continuous. 
\QED

\ber
{\rm (The reason we adopt the $C^\infty$ topology.) 

Lemma \ref{cont} does not hold for the Whitney $C^\infty$ topology. 

For example, in 
$X = C^\infty(\R, \R)$, consider the differentiable mapping 

$h : \R \to C^\infty(\R, \R)$ defined by 

the differentiable mapping $H(\lambda, x) := \lambda$.  
Then $h(0)$ is identically $0$. Its graph is $\R\times 0 \subset 
\R\times\R$. Then there exists an open set $U$ containing 
$\R\times 0$ such that $h^{-1}(W(U)) = \{ 0\}$.

Then $W(U)$ is an open subset of 

$C^\infty(\R, \R)$ with respect to the Whitney $C^\infty$ topology, while 
$h^{-1}(W(U)) = \{ 0\} \subset \R$ is not open in $\R$. 

Therefore 
$h$ is not continuous in the Whitney $C^\infty$ topology. 
}
\enr

\ber
{\rm (The continuity condition in (ii) is necessary.) 

In the above definition (ii), 

the continuity of $k$ is not implied from just the condition that 

for any differentiable mapping 

$h : P \to X$ in the sense (i), the composition 

$k\circ h : P \to Q$ 

is differentiable.

In fact 
set $X = \{ 1/n \} \cup \{ 0 \} \subset \R = C^\infty(\{\pt\}, \R)$ and 
$Y = \{ 0, 1\}  = C^\infty(\{\pt\}, \{ 0, 1\})$.

Define 

$k : X \to Y$ by 
$k(1/n) = 1, k(0) = 0$. 

Then 

any 

differentiable mapping $h : P \to X$ is locally 
constant, and so is $k\circ h : P \to Y$. 

Then $k\circ h$ is differentiable, while 
$k$ is not continuous. 

Thus, in the definition (ii), we need the continuity of $k$. 

}
\enr

\bel
\label{comp}
{\rm (1)} 
The identity 
mapping $\id : X/\!\!\sim\ \to X/\!\!\sim$ is differentiable. 
{\rm (2)}  
Let $\varphi : X/\!\!\sim\  \to Y/\!\!\approx$ 
and $\psi : Y/\!\!\approx \ \to Z/\!\!\equiv$ 
be differentiable mappings. Then 
the composition $\psi\circ\varphi :  X/\!\!\sim\  \to Z/\!\!\equiv$ 
is differentiable. 
\enl

\Proof
(1) is clear since $\id$ is continuous. 
(2) 
Since $\varphi$ and $\psi$ are continuous, $\psi\circ\varphi$ 
is continuous. 
Let $\ell : U (\subseteq Z/\!\!\equiv) \to Q$ 
be a differentiable mapping. 
Then $\ell\circ\psi$ is differentiable by (iv). 
By (iv) again, $(\ell\circ\psi)\circ\varphi = 
\ell\circ(\psi\circ\varphi)$ 
is differentiable. Hence $\psi\circ\varphi$ is differentiable. 
\enP

\bel
{\rm (1)} The quotient mapping $\pi : X \to X/\!\!\sim$ is 
differentiable. 
{\rm (2)} A mapping $\varphi : X/\!\!\sim\   \to Y/\!\!\approx$ 
is differentiable if and only if 
$\varphi\circ\pi : X \to Y/\!\!\approx$ is differentiable. 
\enl

\Proof
(1) First $\pi$ is continuous by the definition of the quotient topology. 
Let $\ell : U (\subseteq X/\!\!\sim) \to Q$ be a differentiable 
mapping. Then $\ell\circ\pi : \pi^{-1}(U) \to Q$ 
is differentiable by the definition (iii). Therefore, by the definition (iv), 
$\pi$ is differentiable. 

(2) First note that $\varphi$ is continuous if and only if 
$\varphi\circ \pi$ is continuous. 
If $\varphi$ is differentiable, then $\varphi\circ\pi$ is differentiable 
by (1) and Lemma \ref{comp}. 
Conversely let $\varphi\circ\pi$ be differentiable, 
and $\ell : U (\subseteq  Y/\!\!\approx) \to Q$ be differentiable. 
Then $\ell\circ(\varphi\circ\pi) = (\ell\circ\varphi)\circ\pi : 
(\ell\circ\varphi\circ\pi)^{-1}(U) \to Q$ 
is differentiable. Therefore $\ell\circ\varphi : 
(\ell\circ\varphi)^{-1}(U) \to Q$ is differentiable. 
Thus, by (iv), $\varphi$ is differentiable. 
\enP

\ber
{\rm
Let $X/\!\!\sim$ be a mapping space quotient 
and $U \ \subseteq X/\!\!\sim$ be an open subset. 
Then 
$$
E(U) := \{ 
\ell : U \to \R \mid \ell {\mbox{\rm \ is differentiable.}} \} 
$$
is an $\R$-algebra. In fact, for $\ell, \ell' \in E(U)$, we have 
$\ell + \ell', \ell\cdot\ell' \in E(U)$. Moreover 
any constant function $U \to \R$ is differentiable. 
Furthermore, for any $\ell_1, \dots, \ell_r  \in E(U)$ 
and for any $C^\infty$ function $\tau : \R^r \to \R$, 
the composition $\tau\circ (\ell_1, \dots, \ell_r)$ belongs to 
$E(U)$. Thus $E(U)$ is a $C^\infty$-ring and 
$E$ is a sheaf of $C^\infty$-rings in the sense of \cite{fourier}. 
}
\enr

We have the following results: 

\bel
For $C^\infty$ manifolds $N, M, L$, the composition 
$c: C^\infty(N, M)\times C^\infty(M, L) 
\to C^\infty(N, L)$ is differentiable. 
\enl

\Proof
First we remark that $c$ is continuous for the $C^\infty$ topology \cite{GG}. 
Let $h : P \to C^\infty(N, M) \times C^\infty(M, L)$ be a differentiable mapping. 
Suppose $H : P\times (N \coprod M) = P\times N \coprod 
P\times M \to (M \coprod L)$ is a $C^\infty$ mapping defining $h$. 
Set $h(p) = (f(p), g(p)), p \in P$.  Then $H(P\times N) \subseteq M$ and 
$H(P\times M) \subseteq L$. 
Then $f(p)(x) = H(p, x)$ for $p \in P, x \in X$, and 
$g(p)(y) = H(p, y)$ for $p \in P, y \in Y$. 
Then $(c\circ h)(p)(x) = g(p)(f(p)(x)) = H(p, H(p, x))$. 
Therefore the $C^\infty$ mapping $K : P\times N \to L$ 
defined by $K(p, x) := H(p, H(p, x))$ defines 
$c\times h$. Therefore $c\circ h$ is differentiable. Therefore, by 
Lemma \ref{diff-equiv}, we see $c$ is differentiable. 
\QED

\bel
If $N$ and $N'$ are diffeomorphic, and 
if $M$ and $M'$ are diffeomorphic, then 
$C^\infty(N, M)$ and $C^\infty(N', M')$ are 
diffeomorphic. 
\enl

\Proof
Let $\sigma : N \to N'$ and 
$\tau : M \to M'$ be 
diffeomorphisms. 
Then set $\varphi : C^\infty(N, M) \to C^\infty(N', M')$ 
by $\varphi(f) := \tau\circ f\circ\sigma^{-1}$. 
Then $\varphi$ is differentiable. 
In fact, let $h : P \to C^\infty(N, M)$ be a differentiable mapping. 
Suppose $h$ is defined by a $C^\infty$ mapping 
$H : P\times N \to M$. 
Then $\varphi\circ h : P \to C^\infty(N', M')$ is defined by the 
$C^\infty$ mapping $\tau\circ H\circ(\id_P\times \sigma^{-1}) 
: P\times N' \to M'$. Therefore $\varphi\circ h$ is differentiable. 
By Lemma \ref{diff-equiv}, we see $\varphi$ is differentiable. 
By symmetry, define $\psi :  C^\infty(N', M') \to C^\infty(N, M)$ 
by $\psi(g) := \tau^{-1}\circ g\circ \sigma$. Then 
$\psi$ is differentiable and it is the inverse of $\varphi$. 
Therefore $\varphi$ is a diffeomorphism. 
\enP

\section{Mild quotients.}
\label{Mild spaces.} 
Apart from general setting, we select \lq\lq mild quotients \rq\rq 
among general mapping space quotients in our purpose. 

A mapping space 
$X \subseteq C^\infty(N, M)$ 
is called a {\it mild space} 
\\
(a): the topology of $X$ has a countable basis, 
and \\ 
(b): 
for any $x_0 \in X$ 
and for any sequence $x_n \in X, (n = 1, 2, 3, \dots)$ with 
$x_n \to x_0 (n \to \infty)$, 
there exist a differentiable mapping 
$h : \R \to X$, a subsequence $x_{n_k}$ 
and a sequence $p_0, p_k \in \R (k = 1, 2,3, \dots)$ 
satisfying $p_k \to p_0 (k \to \infty)$ and 
$h(p_k) = x_{n_k} (k = 1, 2,3, \dots)$. 

A mapping space quotient is called a {\it mild quotient} 
if, it is diffeomorphic to a quotient 
$X/\!\!\sim$ of a mild space $X$. 

\

The motivation of introducing the notion of mild quotient 
lies in the following (cf. Lemma \ref{diff-equiv}): 

\bel
Let $\varphi : X/\!\!\sim \ \to Z/\!\!\equiv$ 
be a mapping between mapping space quotients. 
Suppose that $X/\!\!\sim$ is a mild quotient of a mild space 
$X$ and 
that $\varphi$ satisfying the condition: 
For any differentiable mapping 
$h : P \to X$ 
from a differentiable manifold $P$, 
the composition $\varphi\pi\circ\circ h : 
P \to Z/\!\!\equiv$ is differentiable. 
Then $\varphi$ is a continuous mapping. 
Therefore $\varphi$ is 
differentiable by Lemma \ref{diff-equiv}. 
\enl

\Proof
Assume $\varphi$ is not continuous. 
Then $\varphi\circ\pi : X \to Z/\!\!\equiv$ is not 
continuous. 
Note that 
$\varphi : X/\!\!\sim \ \to Z/\!\!\equiv$ is continuous 
if and only if 
$\varphi\circ\pi : X \to Z/\!\!\equiv$ is continuous. 
We may suppose $X$ itself is a mild space, i.e., $X$ 
satisfies the conditions (a)(b). 
By (a),  $X$ has the topology with countable bases. 

Then there exist 
a point $x_0 \in X$ and a sequence 
$x_n \in X$ with $x_n \to x_0 (n \to \infty)$ 
while the sequence $(\varphi\circ\pi)(x_n)$ does not converge 
to $(\varphi\circ\pi)(x_0) \in Y/\!\!\approx$; 
there exists an open set $V \subset Y/\!\!\approx$ 
with $(\varphi\circ\pi)(x_0) \in V, 
(\varphi\circ\pi)(x_n) \not\in V, (n = 1, 2, 3, \dots)$. 
Since $X$ is a mild space, there exist a 
differential mapping   
$h : \R \to X$, a subsequence $x_{n_k}$ 
and a sequence $p_k \in \R$ with $p_k 
\to p_0 \in P$ and 
$h(p_k) = x_{n_k}$. 
Then $(\varphi\circ\pi)\circ h : 
P \to Y/\!\!\approx$ must be continuous. 
Therefore $((\varphi\circ\pi)\circ h)(p_k) = (\varphi\circ\pi)(x_{n_k})$ 
converges to $((\varphi\circ\pi)\circ h)(p_0) = (\varphi\circ\pi)(h(p_0)) = (\varphi\circ\pi)(x_0)$. 
This leads to a contradiction. 
\enP

\bee
{\rm
A $C^\infty$ manifold is a mild space. 
}
\ene

\bee

\label{mapping space is mild}
{\rm (\cite{Izumi}) 
Let $N, M$ be differentiable manifolds. Then 
$C^\infty(N, M)$ is a mild space. In fact 
let $f_n \in C^\infty(N, M), (n = 1, 2, 3, \dots)$ be a sequence 
converging to $f_0 \in C^\infty(N, M)$ for the $C^\infty$ 
topology. 
Then, using Hestence's lemma \cite{Tougeron}, we see that 
there exist a positive sequence $a_k$, 
a subsequence $f_{n_k}$ of $f_n$, 
a $C^\infty$ mapping $H : \R \times N 
\to M$ satisfying $a_k \to 0$ and  $H(a_k, x) = f_{n_k}(x),
 (x \in N, n = 1, 2, 3, \dots)$. 
Then we have the differentiable mapping $h : \R 
\to C^\infty(N, M)$ defined by 
$h(t)(x) = H(t, x), (t \in \R, x \in X)$ and that 
$h(a_k) = f_{n_k}$ with $a_k \to 0$. 
}
\ene
 
\bee
\label{Diff_c}
{\rm 
Let $M$ be a manifold. Let $\Diff(M)$ (resp. 
$\Diff_c(M)$) 
denotes the space of all diffeomorphisms 
(resp. diffeomorphisms with compact supports). 
Then $\Diff_c(M)$ is a mild space. 
(I do not know whether $\Diff(M)$ is mild or not 
if $M$ is non-compact. ) 
}
\ene

\bee
{\rm
The subspace $\QQ \subset \R$, the set of rational numbers in the real, 
is not mild. 
}
\ene

\bee
\label{Open subsets of a mild space are mild}
{\rm 
Open subsets of a mild quotient are mild: 
Let $\bar{X} = X/\!\!\sim$ be a quotient space of 
of a mild space $X \subset C^\infty(N, M)$ 
by the projection $\pi : X \to \bar{X}$. 
Let $W \subset \bar{X}$ be an open subset 
of  $\bar{X}$. 
Set $X' := \pi^{-1}(W)$. 
Then $W$ is regarded as a quotient space by $\pi : X'
\to W$. 
Then $W$ is a mild space. 
}
\ene

The following is easy to verify:

\bel
\label{Product of mild spaces is mild}
Let $X \subseteq C^\infty(N, M)$ and 
$Y \subseteq C^\infty(L, \Lambda)$ be mild spaces. 
Then the product 
$X \times Y \subseteq 
C^\infty(N\coprod L, M \coprod \Lambda)$ 
is a mild space. 
\enl

\section{Differentiable structure on a moduli space.}
\label{Differentiable structure on a moduli space.}

Recall that we have defined the labelled moduli space of 
$f : S^1 \to \R^2$ as the mapping space quotient 
$$
\widetilde{\mathcal M}_\symp(f) = 
\Diff^+(\R^2)/\!\!\sim_f, 
$$
of $\Diff^+(\R^2) \subset C^\infty(\R^2, \R^2)$. 
Therefore naturally we define the differentiable structure 
on $\widetilde{\mathcal M}_\symp(f)$. 
Also we can consider another mapping space quotient 
$
\Diff_c(\R^2)/\!\!\sim_f
$
of $\Diff_c(\R^2) \subset \Diff^+(\R^2) \subset C^\infty(\R^2, \R^2)$. 
In fact we have the following: 

\bel
\label{compact support}
$\Diff^+(\R^2)/\!\!\sim_f$ is diffeomorphic to 
$\Diff_c(\R^2)/\!\!\sim_f$. 
\enl

\Proof
Let $D_1 = D(R) \subset \R^2$ be a closed disk 
with radius $R$ whose interior contains $f(S^1)$. 
Set $D_2 = D(2R)$. 
We denote by $\Diff_D(\R^2)$ the set of diffeomorphisms on $\R^2$ 
with support contained in $D$. 
Then we show that 
$$
\Diff^+(\R^2)/\!\!\sim_f \ 
\cong
\Diff_{D_2}(\R^2)/\!\!\sim_f \ 
\cong
\Diff_c(\R^2)/\!\!\sim_f. 
$$
Denote by $i : \Diff_{D_2}(\R^2) \to \Diff^+(\R^2)$ 
the inclusion. Then $i$ induces the mapping 
$\varphi : \Diff_{D_2}(\R^2)/\!\!\sim_f \to 
\Diff^+(\R^2)/\!\!\sim_f$ which is differentiable and bijective. 
In fact we define the inverse $\psi : \Diff^+(\R^2)/\!\!\sim_f \to 
\Diff_{D_2}(\R^2)/\!\!\sim_f$ by using 
a differentiable mapping $r : \Diff^+(\R^2) \to \Diff_{D_2}(\R^2)$ 
with $r\vert_{\Diff_{D_1}(\R^2)} = \id$. 
Therefore $\varphi$ is a diffeomorphism. 
A diffeomorphism between $\Diff_c(\R^2)/\!\!\sim_f$ and 
$\Diff_{D_2}(\R^2)/\!\!\sim_f$ is obtained similarly. 
\QED

\ 

Recall that the symplectic moduli space has been defined by 
${\mathcal M}_\symp(f) = \Diff^+(\R^2)f/\Symp(\R^2)$. 
Therefore it is regarded as a quotient of a subspace of 
$C^\infty(S^1, \R^2)$. 
In fact $\Diff^+(\R^2)f \subset C^\infty(S^1, \R^2)/\Diff^+(S^1)$. 
Note that $f \in C^\infty(S^1, \R^2)/\Diff^+(S^1)$. 
On the other hand, since there is a canonical surjection 
$\Pi : \widetilde{\mathcal M}_\symp(f) \to 
{\mathcal M}_\symp(f)$ defined by $\Pi([\rho]) = [\rho\circ f]$, 
the space ${\mathcal M}_\symp(f)$ can be regarded as a quotient 
of $\Diff^+(\R^2) \subset C^\infty(\R^2, \R^2)$ as well. 
However the differentiable structure of ${\mathcal M}_\symp(f)$ 
does not depend on these representations as mapping space quotients as we 
will see below. 

We have set, for an oriented curve $f \in C^\infty(S^1, \R^2)/\Diff^+(S^1)$, 
the symmetry group of $f$, 
$$
S_f := \{ \rho \in \Diff^+(\R^2) \mid \rho\circ f = f \ \mod. \ \Diff^+(S^1) 
\}. 
$$
Moreover we set, the groups of symmetries with compact supports, 
$$
S_{c,f} := \{ \rho \in \Diff_c(\R^2) \mid \rho\circ f = f \ \mod. \ \Diff^+(S^1) 
\} . 
$$
Then, in fact we have the following: 

\bel
If $f$ is of finite type, then there are diffeomorphisms 
$$
\begin{array}{c}
\Diff_c(\R^2)/S_{c,f} \cong \Diff^+(\R^2)/S_{f} 
\cong \Diff^+(\R^2)f 
\vspace{0.2truecm}
\\
\quad\quad\quad (= \Diff_c(\R^2)f 
\subset C^\infty(S^1, \R^2)/\Diff^+(S^1))
\end{array}
$$
\enl

\Proof
The first diffeomorphism is given similarly as in Lemma \ref{compact support}. 
To give the second diffeomorphism 
define $\widetilde{\varphi} : \Diff^+(\R^2)
\to \Diff^+(\R^2)f$ by $\varphi(\rho) = \rho\circ f$. 
Then $\varphi$ is differentiable and 
$\varphi$ induces a differentiable mapping $\varphi : 
\Diff^+(\R^2)/S_f \to \Diff^+(\R^2)f$, which is bijective. 
To see the inverse of $\varphi$ is differentiable, 
we take any differentiable mapping 
$h : P \to C^\infty(S^1, \R^2)$ 
which induces a differentiable mapping $\bar{h} : 
P \to \Diff^+(\R^2)f$, and we show that there exists 
a differentiable mapping $\widetilde{h} : P \to \Diff^+(\R^2)$ 
that $\widetilde{\varphi} \circ \widetilde{h} = \bar{h}$. 
To do that, for a finite dimensional $C^\infty$ family $f_\lambda$ 
isotopic to $f$, $\lambda \in P$,  
we have to find a finite dimensional $C^\infty$ family of 
diffeomorphisms $\rho_\lambda$ on $\R^2$ and $\sigma_\lambda$ 
on $S^1$ such that 
$\rho_\lambda\circ f \circ \sigma_\lambda 
= f_\lambda$. (An isotopy is covered by 
an ambient isotopy). Since $f$ is of finite type, any multi-germ of 
$f$ is finitely determined(\cite{Mather3}\cite{Wall}). 
Therefore, any $C^\infty$ family in an orbit 
is covered by a $C^\infty$ family of diffeomorphisms  
locally at any point on $f(S^1)$. By patching local ambient isotopies, 
we find a global ambient isotopy $(\rho_\lambda, \sigma_\lambda)$. 
\QED

\ 

Thus we have the following basic result: 

\bel
\label{moduli space is mild}
Let $f \in C^\infty(S^1, \R^2)/\Diff^+(S^1)$ be an oriented curve. 
Then 
the labelled symplectic moduli space 
$\widetilde{\mathcal M}_\symp(f) = \Diff^+(\R^2)/\sim_f$ 
is a mild quotient. 
If $f : S^1 \to \R^2$ is of finite type, then 
$\Diff^+(\R^2)f \subseteq C^\infty(S^1, \R^2)/\Diff^+(S^1)$ 
is a mild quotient. Therefore 
%the symplectic moduli space 
${\mathcal M}_\symp(f) = \Diff^+(\R^2)f/\Symp(\R^2)$ 
is a mild quotient. 
\enl

\ber
{\rm
The space of germs $C^\infty((N, S), (M, y_0))$ is 
a mild quotient, since it is a quotient of the mild space 
$C^\infty(N, M)$. 
Similarly to Lemma \ref{moduli space is mild}, we see 
also the (labelled) local symplectic 
moduli space of a map-germ of finite codimension is 
a mild quotient. 
}
\enr

\section{The space of symplectic forms.}
\label{The space of symplectic forms.}

Differential two-forms $\omega, \omega'$ on $\R^2$ are equivalent 
if there exists a diffeomorphism $\rho$ on $\R^2$ such that 
$\omega = \rho^*\omega'$. 
Let ${\mathcal S}(\R^2)$ denote 
the space of symplectic forms on $\R^2$ which are equivalent to 
the standard symplectic form $\omega_0 = dx\wedge dy$. 
It is known that $\omega \in {\mathcal S}(\R^2)$ 
if and only if 
$\displaystyle{\int_{\R^2} \omega = \pm \infty}$ (cf. \cite{GS}). 
Note that ${\mathcal S}(\R^2)$ is a subset 
of $C^\infty(\R^2, \bigwedge^2(T^*\R^2))$. 
A symplectic form on $\R^2$ is called {\it 
positive} 
if it gives the standard orientation on $\R^2$. 
Then also consider the subspace 
${\mathcal S}^+(\R^2) \subset {\mathcal S}(\R^2)$ 
consisting of positive forms 
on $\R^2$. Note that the standard symplectic 
form $\omega_0 = dx \wedge dy$ belongs to 
${\mathcal S}^+(\R^2)$. 

Actually we treat another object ${\mathcal S}_c(\R^2)$, 
the space of symplectic forms on $\R^2$ with 
\lq\lq compact support\rq\rq, 
namely symplectic forms which agree with 
$\omega_0$ outside of compact subsets. 
Note that ${\mathcal S}_c(\R^2) \subset {\mathcal S}^+(\R^2)$. 
Moreover we see that, for any $\omega \in {\mathcal S}_c(\R^2)$,  
there exists a $\rho \in \Diff(\R^2)$ with compact support 
satisfying $\omega = \rho^*\omega_0$. 
We consider the group $\Diff_c(\R^2)$ (resp. 
$\Symp_c(\R^2)$) 
of diffeomorphisms (resp. symplectomorphisms) 
on $\R^2$ with compact supports. 
We define an equivalence relation $\sim$ on $\Diff_c(\R^2)$,  
using the subgroup $\Symp_c(\R^2)$; $\rho_1, \rho_2 
\in \Diff_c(\R^2)$ are equivalent, $\rho_1 \sim \rho_2$, 
if $\rho_1 = \tau \circ \rho_2$ 
for some $\tau \in \Symp_c(\R^2)$. 
Then we have the quotient space $\Diff_c(\R^2)/\sim$. 
Though the action is a left action, 
we denote the quotient space by $\Diff_c(\R^2)/\Symp_c(\R^2)$. 
Recall that the notion of mild spaces introduced in 
\S \ref{Mild spaces.}. 
Then, similarly to Example \ref{Diff_c},  we have: 

\bel
${\mathcal S}_c(\R^2) \subset C^\infty(\R^2, \bigwedge^2(T^*\R^2))$
is a mild space. 
Therefore any quotient space of ${\mathcal S}_c(\R^2)$ is a mild quotient. 
\enl

Moreover we have: 

\bel
\label{forms and diffeomorphisms}
${\mathcal S}_c(\R^2)$ is diffeomorphic to 
$\Diff_c(\R^2)/\Symp_c(\R^2)$. 
\enl

\Proof
Define a mapping $\Phi : \Diff_c(\R^2) \to {\mathcal S}_c(\R^2)$ 
by $\Phi(\rho) = \rho^*\omega_0$. Then $\Phi$ is a surjective 
mapping. 
Note that if $\rho_1 = \tau\circ\rho_2$ 
with $\tau \in \Symp_c(\R^2)$, then 
$\rho_1^*\omega_0 = \rho_2^*\omega_0$, therefore 
$\rho_1\circ\rho_2^{-1} \in \Symp_c(\R^2)$. 
The mapping $\Phi$ induces a mapping 
$\varphi : \Diff_c(\R^2)/\Symp_c(\R^2) \to 
{\mathcal S}_c(\R^2)$ 
defined by $\varphi([\rho]) = \rho^*\omega_0$. 
Then $\varphi$ is a bijection. 
Let $h : P \to \Diff_c(\R^2)$ be any differentiable mapping from 
any manifold. Let $H : P\times \R^2 \to \R^2$ be a $C^\infty$ mapping which defines $h$. Then we set 
$H' : P \times \R^2 \to \bigwedge^2(T^*\R^2)$ 
by $H'(p, x, y) = (h(p)^*\omega)(x, y)$. Then $H'$ is $C^\infty$ and defines 
$\Phi\circ h$. Therefore $\Phi\circ h$ is differentiable. 
Since $\Diff_c(\R^2)$ is a mild space, we see $\Phi$ is differentiable 
by Lemma \ref{diff-equiv}. Therefore $\phi$ is differentiable. 

Conversely define 
$\widetilde{\psi} : 
{\mathcal S}_c(\R^2) \to \Diff_c(\R^2)$ by 
$$
\displaystyle{\widetilde{\psi}(\omega)(x, y) := (\int_0^x f(x, y)dx, \ y)}
$$
where $\omega = f(x, y)dx\wedge dy$. 
Then $\widetilde{\psi}$ induces a mapping 
$$
\psi : 
{\mathcal S}_c(\R^2) \to \Diff_c(\R^2)/\Symp_c(\R^2). 
$$
We see that $\psi$ is the inverse of $\varphi$. 
To show $\psi$ is differentiable, 
let 
$k : P \to {\mathcal S}_c(\R^2)$ be a differentiable mapping. 
Let $K : P\times \R^2 \to \bigwedge^2(T^*\R^2)$ 
be a $C^\infty$ mapping which defines $k$. 
Set 
$$
\displaystyle{K'(p, x, y) = (\int_0^x f(p, x, y)dx, \ y)}
$$
where $K'(p, x, y) = f(p, x, y)dx\wedge dy$. 
Then $K' : P\times \R^2 \to \R^2$ is $C^\infty$ and $K'$ defines 
$\widetilde{\psi}\circ k : P \to \Diff_c(\R^2)$. 
Therefore 
$\widetilde{\psi}\circ k$ is differentiable. 
Thus $\psi\circ k = \pi\circ(\widetilde{\psi}\circ k) : 
P \to \Diff_c(\R^2)/\Symp_c(\R^2)$ 
is differentiable. 
Since ${\mathcal S}_c(\R^2)$ is a mild space, 
we see $\psi$ is continuous. Therefore $\psi$ is differentiable, 
by Lemma \ref{diff-equiv}. 
Thus we show that $\varphi$ is a diffeomorphism. 
\enP

For a oriented plane 
curve $f \in C^\infty(S^1, \R^2)/\Diff^+(S^1)$, 
we have considered the {\it group of symmetry} of $f$:  
$$
S_f := \{ \rho \in \Diff^+(\R^2) \mid \rho\circ f = f  
\ \ {\mbox{\rm up to }} \Diff^+(S^1) \}. 
$$ 
Furthermore, if $f$ is labelled, then 
we consider the {\it group of label-preserving symmetry} of $f$: 
$$
\underline{S}_f := 
\{ \rho \in S_f \mid \rho 
{\mbox{\rm 
\  preserves the labelling of }} f. \}. 
$$
Moreover set $S_{c, f} := S_f \cap \Diff_c(\R^2)$ 
and $\underline{S}_{c, f} := \underline{S}_f \cap \Diff_c(\R^2)$.

Note that there exist exact sequences of groups: 
$$
1 \rightarrow \underline{S}_f \rightarrow S_f \rightarrow G_f \rightarrow 1, 
\quad 
1 \rightarrow \underline{S}_{c, f} \rightarrow S_{c, f} \rightarrow G_f \rightarrow 1. 
$$

On the other hand, let $\omega \in {\mathcal S}^+(\R^2)$. 
Then there exists $\rho \in \Diff^+(\R^2)$ such that 
$\omega = \rho^*\omega_0$. If $\rho_1^*\omega_0 = \rho_2^*\omega_0$, 
then $\tau := \rho_1\circ\rho_2^{-1} \in  \Symp(\R^2)$, and 
$\tau\circ \rho_2 = \rho_1$, therefore, $\tau\circ \rho_2\circ f = \rho_1\circ f$. 
Thus a mapping  
$$
p : {\mathcal S}^+(\R^2) \to \widetilde{\mathcal M}_\symp(f)
$$ 
is well-defined 
by $p(\omega) := [\rho] \in \Diff^+(\R^2)/\!\!\sim_f$ for some $\rho \in \Diff^+(\R^2)$ with 
$\omega = \rho^*\omega_0$. 
Moreover we have the diagram: 
$$
\begin{array}{ccccccc}
\Diff_c(\R^2) & \stackrel{i}{\rightarrow} & \Diff^+(\R^2) & 
\stackrel{\Pi}{\rightarrow} & 
\widetilde{\mathcal M}_\symp(f) & \stackrel{\pi}{\rightarrow} & {\mathcal M}_\symp(f) 
\vspace{0.2truecm}
\\
\Phi\downarrow & & \Phi\downarrow & \ \ \nearrow p & & & 
\vspace{0.2truecm}
\\
{\mathcal S}_c(\R^2) & \stackrel{j}{\rightarrow} & {\mathcal S}^+(\R^2) & 
& & & 
\end{array}
$$

Here $i$ and $j$ are inclusions, 
$\Pi$ and $\pi$ are projections. 
We set $\Phi(\rho) = \rho^*\omega_0$, for $\omega \in \Diff^+(\R^2)$. 

Then we have 

\bel
\label{forms and diff}

{\rm (1)} There are diffeomorphisms 
$$
\begin{array}{c}
{\mathcal S}^+(\R^2)/S_{f}
\cong 
{\mathcal S}_c(\R^2)/S_{c, f}
\cong 
\Diff_c(\R^2)f/\Symp_c(\R^2) 
\vspace{0.2truecm}
\\ 
\cong 
\Diff^+(\R^2)f/\Symp(\R^2) 
=: 
{\mathcal M}_\symp(f). 
\end{array}
$$

{\rm (2)} There are diffeomorphisms 

$$
{\mathcal S}^+(\R^2)/\underline{S}_{f} 
\cong
{\mathcal S}_c(\R^2)/\underline{S}_{c, f} 
\cong 
\Diff_c(\R^2)/\!\!\sim_f \ 
\cong 
\Diff^+(\R^2)/\!\!\sim_f \ 
=: \widetilde{\mathcal M}_\symp(f). 
$$
\enl

\Proof
For $\omega, \omega' \in {\mathcal S}_c(\R^2)$, 
suppose $p\circ j(\omega_1) = p\circ j(\omega_2)$. Then 
$\omega_1 = \rho_1^*\omega_0, \omega_2 = \rho_2^*\omega_0$ 
and $\tau\circ \rho_1\circ f = \rho_2\circ f$ for some 
$\rho_1, \rho_2 \in \Diff_c(\R^2), \tau \in \Symp_c(\R^2)$, 
$\tau$ being label-preserving. 
Then $\rho_2^{-1}\circ\tau \circ\rho_1 \in \underline{S}_{c, f}$ and 
$(\rho_2^{-1}\circ\tau \circ\rho_1)^*\omega_2 = \omega_1$. 
Conversely suppose $\rho^*\omega_2 = \omega_1$, 
for some $\rho \in \underline{S}_{c, f}$. Then 
$\rho_2\circ\rho\circ\rho_1^{-1} =: \tau \in \Symp_c(\R^2)$. 
Since $\rho = \rho_2^{-1}\circ\tau\circ\rho_1 \in \underline{S}_{c, f}$, 
we have $\tau\circ\rho_1\circ f = \rho_2\circ f$, 
$\tau \in \Symp_c(\R^2)$ and $\tau$ preserves the labelling. 
Thus $p\circ j$ induces a bijection $\overline{p} : 
{\mathcal S}_c(\R^2)/\underline{S}_{c, f} \to \Diff_c(\R^2)/\!\!\sim_f$. 
Then $\overline{p}$ is differentiable. In fact, 
by Lemma \ref{forms and diffeomorphisms}, $\overline{p}$ 
is induced from the differentiable mapping $\Pi\circ i$. 
Moreover 
$\psi : \Diff_c(\R^2)/\!\!\sim_f \to {\mathcal S}_c(\R^2)/\underline{S}_{c, f}$ 
defined by $\psi([\rho]) := [\rho^*\omega_0]$ is differentiable and 
$\psi$ is the inverse of $\overline{p}$. 
Therefore $\overline{p}$ is a diffeomorphism. 
For other diffeomorphisms are obtained similarly. 
\QED

\

Let us denote by ${\mathcal S}^+(\R^2, y_0)$ the space of germs of 
positive symplectic forms on $(\R^2, y_0)$.  
Then similarly to Lemma \ref{forms and diff}, we have, 
on local moduli spaces: 

\bel
{\rm (1)}
For a map-germ $f_{y_0} : (S^1, S) \to (\R^2, y_0)$ 
of finite codimension,  
$$
{\mathcal S}^+(\R^2, y_0)/S_{f_{y_0}}
\cong 
\Diff^+(\R^2, y_0)f/\Symp(\R^2, y_0) 
\cong 
{\mathcal M}_\symp(f_{y_0}). 
$$

{\rm (2)} 
For a labelled map-germ $f_{y_0} : (S^1, S) \to (\R^2, y_0)$ 
of finite codimension, 
$$
{\mathcal S}^+(\R^2, y_0)/\underline{S}_{f_{y_0}} 
\cong 
\Diff^+(\R^2, y_0)/\!\!\sim_{f_{y_0}} \ 
=: \widetilde{\mathcal M}_\symp(f_{y_0}). 
$$
\enl

\

\noindent {\it Proof of Theorem \ref{generic immersion}:} 
Let 
$f_0 \in C^\infty(S^1, \R^2)/\Diff^+(S^1)$ be 
a generic immersion. Define 
$$
\Phi :  {\mathcal S}_c(\R^2) \to H^2(\R^2, f_0(S^1), \R) 
= H_2(\R^2, f_0(S^1), \R)^*
$$
by setting, for $\omega \in {\mathcal S}_c(\R^2)$, 
$$
\Phi(\omega)([D_j]) := \int_{D_j} \omega, 
$$
$([D_1], \dots, [D_r] \in H_2(\R^2, f_0(S^1), \R))$. 
Then, via the diffeomorphism 
$$
{\mathcal S}_c(\R^2)/S_{c,f_0} \cong 
\widetilde{\mathcal M}_\symp(f),
$$
the mapping 
$\Phi$ induces the mapping $\varphi$ introduced in 
Introduction. In fact, if $\omega = \rho^*\omega_0$, 
then 
$$
\int_{D_j} \omega = \int_{D_j} \rho^*\omega_0 = 
\int_{\rho(D_j)} \omega_0. 
$$

Clearly the 
image of $\Phi$ is contained in the positive cone 
$H^2(\R^2, f_0(S^1), \R)_{>0}$. 
Moreover, we show that \\
(a): 
if $\alpha$ is in the image of $\Phi$, and $c > 0$, 
then $c\alpha$ is also in the image of $\Phi$, and 
\\
(b): 
if $\alpha$ is in the image of $\Phi$, 
and $\beta \in H^2(\R^2, f_0(S^1), \R)_{> 0}$, then 
$\alpha + \beta$ belongs to the image of $\Phi$. 

In fact, if $\Phi(\omega) = \alpha$, then $\Phi(c\omega) = c\alpha$. 
So we have (a). 

To show (b), 
let $\alpha = \Phi(\omega) \in H^2(\R^2, f_0(S^1), \R)_{> 0}$. 
Let $D$ be one of bounded domains surrounded by $f_0(S^1)$. 
Let $(x_0, y_0)$ be an interior point of $D$. Let $\varepsilon > 0$ 
satisfy $D_{\varepsilon}(x_0, y_0) \subset D$, where 
$D_{\varepsilon}(x_0, y_0)$ means the $\varepsilon$ closed disc 
centred at $(x_0, y_0)$. 
Let $\lambda_{(x_0, y_0, \varepsilon)}$ be a non-negative 
$C^\infty$ 
function on $\R^2$ with 
$\supp(\lambda_{(x_0, y_0, \varepsilon)}) 
= D_{\varepsilon}(x_0, y_0)$ and satisfying 
$$
\int_{D}\lambda_{(x_0, y_0, \varepsilon)}\omega_0 = 1. 
$$
Set $\lambda_D := \lambda_{(x_0, y_0, \varepsilon)}$. 
Then, for each $\beta  
\in H^2(\R^2, f_0(S^1), \R)_{> 0}$, set 
$$
\omega_{\alpha + \beta} := \omega + 
\sum_{j=1}^r \beta([D_j])\lambda_{D_j} \omega. 
$$
Then $\omega_{\alpha + \beta} \in {\mathcal S}_c(\R^2)$. 
Moreover we have
$$
\begin{array}{rcl}
\displaystyle{\int_{D_k}\omega_{\alpha+\beta}} & = & 
\displaystyle{\int_{D_k}\omega + \sum_{j=1}^r \beta(D_j)\int_{D_k}\lambda_{D_j}\omega} 
\vspace{0.1truecm}
\\
& = & \alpha([D_k]) + \beta([D_k]) 
\  = \ (\alpha + \beta)([D_k]). 
\end{array}
$$
Therefore we have $\Phi(\omega_{\alpha + \beta}) = \alpha + \beta$. 
Hence we have (b). 

By (a)(b), we 
see that the image of $\Phi$ coincides with 
$H^2(\R^2, f_0(S^1), \R)_{>0}$. 

If $\rho \in \underline{S}_{c, f_0}$, then 
$\Phi(\rho^*\omega) = \Phi(\omega)$. 
Thus $\Phi$ induces 
the surjective mapping 
$\varphi : {\mathcal S}_c(\R^2)/\underline{S}_{c, f_0} 
\to H^2(\R^2, f_0(S^1), \R)_{>0}$.

We will show $\Phi$ is differentiable. To see this, we 
let $h : P \to {\mathcal S}_c(\R^2)$ 
be a differentiable mapping. Let a 
$C^\infty$ mapping 
$H : P\times \R^2 \to 
\bigwedge^2(T^*\R^2)$ define $h$. 
Then 
$\Phi\circ h : P \to H^2(\R^2, f_0(S^1), \R)_{> 0}$ 
is defined by 
$$
(\Phi\circ h)(p)([D_j]) := \int_{D_j} h(p)^*\omega_0. 
$$
Therefore $\Phi\circ h$ is differentiable. Since 
${\mathcal S}_c(\R^2)$ is mild, we see $\Phi$ is differentiable. 
Therefore $\varphi$ is differentiable. 

We are going to show the mapping  
$$
\varphi : {\mathcal S}_c(\R^2)/\underline{S}_{c, f_0} \to 
H^2(\R^2, f_0(S^1), \R)_{> 0}
$$ is actually a diffeomorphism. 
We need to show the surjection 
$\varphi$ is an injection and its inverse 
%$$
%\psi = \varphi^{-1} : H^2(\R^2, f_0(S^1), \R)_{> 0} \to 
%{\mathcal S}_c(\R^2)/\underline{S}_{c, f_0}
%$$ 
is a differentiable mapping.

Suppose $\Phi(\omega) = \Phi(\omega')$, for $\omega, \omega' 
\in {\mathcal S}_c(\R^2)$. 
Set $\omega = f(x, y)dx \wedge dy$ and 
$\omega' = g(x, y) dx \wedge dy$. Here 
$f$ and $g$ are positive functions on $\R^2$ 
which agree with $1$ outside a compact subset. 
Then set 
$f_t := (1 - t)f + tg, 0 \leq t \leq 1$, and 
$\omega_t := f_t(x, y) dx \wedge dy$. 
Then, by Moser's theorem (\cite{Moser}), 
we can find a $C^\infty$ family $\rho_t \in \underline{S}_{c, f_0}$ 
such that $\omega_t = \rho_t^*\omega_0$. 
Therefore $\omega' = \omega_1 = \rho_1^*\omega_0 = 
(\rho_1\circ\rho_0^{-1})^*\omega$. 
Thus $\varphi$ is an injection. 
Therefore $\varphi$ 
is a bijection. 
Let 
$$
\psi := \varphi^{-1} : 
H^2(\R^2, f_0(S^1), \R)_{> 0} \to 
{\mathcal S}_c(\R^2)/\underline{S}_{c, f_0}
$$
be the inverse of $\varphi$. 
We set $\Psi(\alpha + \beta) = \omega_{\alpha + \beta}$ 
as defined above. 
Then $\Psi$ is a differentiable map and it gives a local 
differentiable lifting of $\psi$ on an open set 
$\{ \alpha + \beta \mid \beta > 0\}$ 
for any $\alpha \in H^2(\R^2, f_0(S^1), \R)_{> 0}$. 
Therefore $\psi$ is differentiable. This shows that 
$\varphi$ is a diffeomorphism. 
\enP

\noindent {\it Proof of Theorem \ref{of finite type}:} 
Now let $f \in C^\infty(S^1, \R^2)/\Diff^+(S^1)$ be of finite type. 
We set 
$$
M := \left( {\prod}_{y_0 \in f(S^1)} 
{\mathcal S}^+(\R^2, y_0)/\underline{S}_{f_{y_0}} \right) 
\times H^2(\R^2, f(S^1), \R)_{> 0}. 
$$
Define the mapping 
$$
\widetilde{\Phi} : 
{\mathcal S}_c(\R^2) \to M, 
$$
by
$$
\widetilde{\Phi}(\omega) := (([\omega]_{y_0})_{y_0 \in f(S^1)}, \ 
\varphi([\omega])), 
$$
$$
\varphi([\omega])([D_j]) = \int_{D_j} \omega, \ (1 \leq j \leq r). 
$$
Then we see $\widetilde{\Phi}$ is differentiable. 
The mapping 
$\Phi$ defined in Introduction is given, in term of 
symplectic forms, by 
$$
\Phi : {\mathcal S}_c(\R^2)/\underline{S}_{c, f} \to M, \quad  
\Phi([\omega]) = \widetilde{\Phi}(\omega).
$$
Note that  
both ${\mathcal S}_c(\R^2)/\underline{S}_{c, f}$ 
and $M$ are mild quotients (see \S \ref{Mild spaces.}). 

We are going to show that 
$\Phi$ is a diffeomorphism. 
Since $\widetilde{\Phi}$ is differentiable, 
we see that $\Phi$ is differentiable. 

To see $\Phi$ is injective, 
suppose that $\Phi([\omega_1]) = \Phi([\omega_2])$, 
for $\omega_1, \omega_2 \in {\mathcal S}_c(\R^2)$. 
Then, for any $y_0 \in f(S^1)$, 
$[\omega_1]_{y_0} = [\omega_2]_{y_0} \in 
{\mathcal S}^+(\R^2, y_0)/\underline{S}_{f_{y_0}}$, 
and  $\int_{D_j} \omega_1 = \int_{D_j} \omega_2$, 
$1 \leq j \leq r$. 
We need to take care of unstable points of $f$, which are a finite number of 
points by the finiteness of $f$. 
Set 
$$
\Uns(f) := \{ y_0 \in f(S^1) \mid f_{y_0} \ {\mbox{\rm is unstable}} \}.
$$ 
Then there exist germs of label preserving symmetries 
$\rho_{y_0} : (\R^2, y_0) \to (\R^2, y_0)$ 
of $f_{y_0} : (S^1, S) \to (\R^2, y_0)$ such that 
$\rho_{y_0}^*\omega_2 = \omega_1$ near $y_0$. 
Note that there exist germs $\sigma_{y_0} \in \Diff^+_0(S^1, f^{-1}(y_0))$
 such that $\rho_{y_0}\circ f = f \circ \sigma_{y_0}$ near $f^{-1}(y_0)$. 
The local symmetries $\rho_{y_0}$, $y_0 \in \Uns(f)$, 
is covered by a global symmetry $\rho$ of $f$, via 
a $\sigma \in \Diff^+(S^1)$ giving $\sigma_{y_0}$ near $f^{-1}(y_0)$. 
Then $\rho^*\omega_2 = \omega_1$ on a neighbourhood of 
$\Uns(f)$. This construction is extended on a neighborhood of $f(S^1)$. 
We extend $\rho$ to a $\rho \in \underline{S}_{c,f}$. 
Now 
$$
\int_{D_j} \rho^*\omega_2 
= \int_{\rho(D_j)} \omega_2 = \int_{D_j} \omega_2 = \int_{D_j} \omega_1. 
$$
Then there exist a diffeomorphism $\rho' \in \underline{S}_{c,f}$ 
such that $\rho'$ is identity on a neighborhood of $f(S^1)$ and 
$\rho'^*(\rho^*\omega_2) = \omega_1$. 
Now $\rho\circ\rho' \in \underline{S}_{c,f}$, so 
we have $[\omega_1] = [\omega_2]$. 
Therefore we see that $\Phi$ is an injection. 

Next we show $\Phi$ is surjective. 
Take a germ of positive form $\omega_{y_0} 
= f_{y_0}(x, y)dx\wedge dy$ at every $y_0 \in \Uns(f)$. 
Then, for any given $\varepsilon_j > 0, 1 \leq j \leq r$, 
$(\omega_{y_0})_{y_0 \in \Uns(f)}$ is extended to 
$\omega \in {\mathcal S}_c(\R^2)$ such that 
$\int_{D_j}\omega < \varepsilon_j, 1 \leq j \leq r$. 
This is established by setting 
$$
\omega = \left[ \ \sum_{y_0 \in \Uns(f)} \mu_{y_0}(x, y)f_{y_0}(x, y) 
+ \varepsilon(x, y) \ \right] dx\wedge dy, 
$$
using a non-negative $C^\infty$ functions $\mu_{y_0}(x, y)$ 
and $\varepsilon$ 
on $\R^2$ for $y_0 \in \Uns(f)$ satisfying 
the following conditions: 
(i) $\mu_{y_0}(x, y)$ is equal to $1$ in a neighbourhood 
of $y_0$ and 
$\varepsilon$ is equal to $0$ 
in a smaller neighborhood of the finite set $\Uns(f)$.  
(ii) $\mu_{y_0}(x, y)$ is equal to $0$ in a neighbourhood of $y_0$ and 
$\varepsilon$ is sufficiently small on each $D_j$, so that 
$\int_{D_j} \omega < \varepsilon_j$. 
(iii) $\varepsilon$ is equal to $1$ outside of 
a compact set. 

\begin{figure}[htbp]
  \begin{center}
 \includegraphics[width=6truecm, height=3truecm, clip]{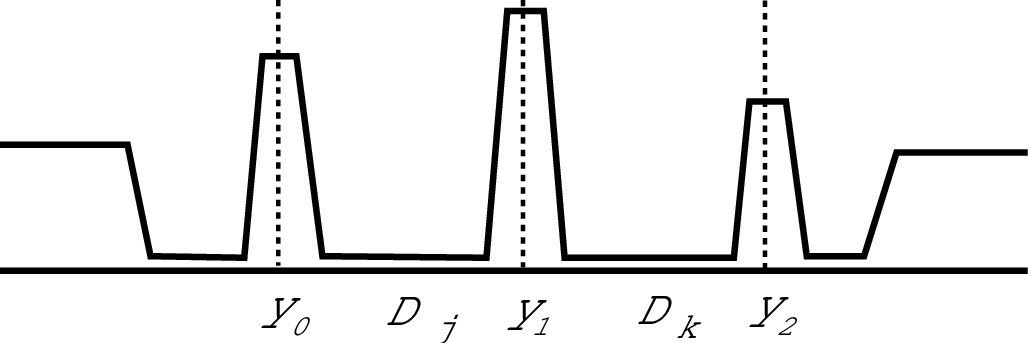} 
 %\vspace{-1.0truecm}
 \caption{Extending form-germs to the plane.}
    \label{bump}
  \end{center}
 \end{figure}

In Fig. \ref{bump}, we illustrate a 
required extension $\omega$ of given form-germs $\omega_{y_0}$. 
Then, for any $\alpha \in H^2(\R^2, f(S^1), \R)_{> 0}$ 
with $\alpha(D_j) > \varepsilon_j$, set 
$$
\omega_{((\omega_{y_0}), \alpha)} = 
\omega + \left( \alpha(D_j) - \int_{D_j}\omega \right) \lambda_{D_j}dx\wedge dy.
$$
Here $\lambda_{D_j}$ is a bump function as introduced in the proof of 
Theorem \ref{generic immersion}. 
Then $\omega_{((\omega_{y_0}), \alpha)} \in {\mathcal S}_c(\R^2)$ 
and 
$$
\widetilde{\Phi}(\omega_{(\omega_{y_0}), \alpha}) 
= ((\omega_{y_0}), \alpha). 
$$
Since $\varepsilon_j > 0$ is arbitrary, we see $\widetilde{\Phi}$ is surjective. 
Thus we see $\Phi$ is bijective.

Let 
$\Psi : M \to {\mathcal S}_c(\R^2)/\underline{S}_{c, f}$ 
be the inverse of $\Phi$. We are going to show $\Psi$ is 
differentiable. 
Let 
$$
h : P \to \left( {\prod}_{y_0 \in \Uns(f)} 
{\mathcal S}^+(\R^2) \right) \times 
H^2(\R^2, f(S^1), \R)_{> 0}
$$
be a differentiable mapping from a finite dimensional manifold $P$. 
We ask whether 
$\Psi\circ h : P \to {\mathcal S}_c(\R^2)/\underline{S}_{c, f}$ 
is differentiable. 
Let a $C^\infty$ mapping  
$$
{\textstyle 
H : P\times (\coprod_{y_0 \in \Uns(f)} \R^2) \coprod \{ \pt \} 
\to (\coprod_{y_0 \in \Uns(f)} \R) \coprod H^2(\R^2, f(S^1), \R)_{> 0}
}
$$
defines $h$, where $\coprod_{y_0 \in \Uns(f)} \R^2$ 
means the disjoint union of $r''$-copies of $\R^2$, 
$r''$ being the number of $\Uns(f)$.  
We are identifying $\R$ with $\bigwedge^2 (\R^2)^*$. 
Set 
$$
H(\lambda, y_1, \dots, y_{r''}, \pt) 
= ((\omega_{y_0, \lambda} 
= f_{y_0}(x, y, \lambda)dx\wedge dy)_{y_0 \in \Uns(f)}, \ \ \alpha_\lambda). 
$$
For each compact subset $\Lambda$ of $P$, we construct 
a $C^\infty$ family $\omega_\lambda \in {\mathcal S}_c(\R^2)$ 
by 
$$
\omega_{((\omega_{y_0, \lambda}), \alpha_\lambda)} := 
\omega + \sum_{j=1}^r \left( \alpha(D_j) - \int_{D_j}\omega \right)
\lambda_{D_j}dx\wedge dy, 
$$
where 
$$
\omega := \left[ \ \sum_{y_0 \in \Uns(f)} \mu_{y_0}(x, y)f_{y_0}(x, y, \lambda) 
+ \varepsilon(x, y) \ \right] dx\wedge dy
$$
is defined by functions $\mu_{y_0}$ and $\varepsilon(x, y)$ 
which are independent of $\lambda \in \Lambda$. 
Then 
$$
\widetilde{\Psi\circ h}({(\omega_{y_0}, \lambda), \alpha_\lambda})
:= \omega_{((\omega_{y_0, \lambda}), \alpha_\lambda)}
$$
gives a local differentiable lifting of $\Psi\circ h$. Therefore 
$\Psi\circ h$ is also differentiable. Since $M$ is a mild quotient, 
we see $\Psi$ is differentiable. 
Thus we have shown that 
$\Phi$ is a diffeomorphism as required. 
\enP

\noindent{\it Proof of Theorems \ref{multi-components} and 
\ref{surface}:} 
The proofs are established in the parallel way to that of Theorem \ref{of finite type}. Therefore they are left to the reader. 
\QED

\begin{flushleft}
Goo ISHIKAWA \\ 
Department of Mathematics, Hokkaido University, 
Sapporo 060-0810, Japan. 
\begin{verbatim}
E-mail : ishikawa@math.sci.hokudai.ac.jp
Web: http://www.math.sci.hokudai.ac.jp/~ishikawa
\end{verbatim} 
\end{flushleft}

\end{document}